\documentstyle[12pt]{article}
\newtheorem{theorem}{{\sc Theorem}}
\newcommand{\bt}{\begin{theorem}}
\newcommand{\et}{\end{theorem}}
\newcommand{\newsection}[1]{\setcounter{equation}{0} \setcounter{theorem}{0}
\section{#1}}

\newcommand{\NI}{\noindent}

\newcommand{\bea}{\begin{eqnarray}}
\newcommand{\eea}{\end{eqnarray}}

\def \b #1 {\bf #1}
\newcommand{\IR}{I\!\!R}
\newcommand{\IE}{I\!\!E}
\newcommand{\IC}{I\!\!C}

\newcommand{\IT}{I\!\!T}
\newcommand{\IN}{I\!\!N}
\newcommand{\IZ}{Z\!\!\!Z}
\newcommand{\clk}{{\cal K}}

\newcommand{\cla}{{\cal A}}

\newcommand{\cli}{{\cal I}}

\newcommand{\clf}{{\cal F}}
\newcommand{\clg}{{\cal G}}
\newcommand{\clh}{{\cal H}}

\newcommand{\clo}{{\cal O}}
\newcommand{\clb}{{\cal B}}

\newcommand{\cln}{{\cal N}}

\newcommand{\clm}{{\cal M}}

\newcommand{\raro}{\rightarrow}

\newcommand{\vsp}{\vskip 1em}

\def \qed {\hfill \vrule height6pt width 6pt depth 0pt}
\newcommand{\be}{\begin{equation}}
\newcommand{\ee}{\end{equation}}
\newcommand{\ben}{\begin{eqnarray*}}
\newcommand{\een}{\end{eqnarray*}}
\begin{document}
\thispagestyle {empty}
\sloppy

\centerline{\large \bf Pure inductive limit state and Kolmogorov's property }

\bigskip
\centerline{\bf Anilesh Mohari }
\smallskip
\centerline{\bf S.N.Bose Center for Basic Sciences, }
\centerline{\bf JD Block, Sector-3, Calcutta-98 }
\centerline{\bf E-mail:anilesh@boson.bose.res.in}
\smallskip
\centerline{\bf Abstract}
\bigskip

\vsp
Let $(\clb,\lambda_t,\psi)$ be a $C^*$-dynamical system where $(\lambda_t: t \in \IT_+)$ be a semigroup of
injective endomorphism and $\psi$ be an $(\lambda_t)$ invariant state on the $C^*$ subalgebra $\clb$ and $\IT_+$ is 
either non-negative integers or real numbers. The central aim of this exposition is to find a useful criteria for 
the inductive limit state $\clb \raro^{\lambda_t} \clb$ canonically associated with $\psi$ to be pure. We achieve this 
by exploring the minimal weak forward and backward Markov processes associated with the Markov semigroup on the corner 
von-Neumann algebra of the support projection of the state $\psi$ to prove that Kolmogorov's property [Mo2] of 
the Markov semigroup is a sufficient condition for the inductive state to be pure. As an application of this criteria
we find a sufficient condition for a translation invariant factor state on a one dimensional quantum spin chain 
to be pure. This criteria in a sense complements criteria obtained in [BJKW,Mo2] as we could go beyond lattice 
symmetric states.

\newpage
\newsection{ Introduction:}

\vsp
Let $\tau=(\tau_t,\;t \ge 0)$ be a semigroup of identity preserving completely
positive maps [Da1,Da2,BR] on a von-Neumann algebra $\cla_0$ acting on a Hilbert space $\clh_0$, 
where either the parameter $t \in \!R_+$, the set of positive real numbers or $ t \in \!Z^+$, the set 
of positive integers. We assume further that the map $\tau_t$ is normal for each $t \ge 0$ and the map 
$t \raro \tau_t(x)$ is weak$^*$ continuous for each $x \in \cla_0$.

\vsp
We say a projection $p \in \cla_0$ is {\it sub-harmonic and harmonic  } if $\tau_t(p) \ge p$ and $\tau_t(p)=p$
for all $t \ge 0$ respectively. For a sub-harmonic projection $p$, we define the reduced quantum dynamical 
semigroup $(\tau^p_t)$ on the von-Neumann algebra $p\cla_0p$ by $\tau^p_t(x)=p\tau_t(x)p$ where $t \ge 0 $ and 
$x \in \cla^p_0$. $1$ is an upper bound for the increasing positive operators $\tau_t(p),\; t \ge 0$. Thus there 
exists an operator $0 \le y \le 1$ so that $y=\mbox{s.lim}_{t \raro \infty} \tau_t(p)$. A normal state $\phi_0$ 
is called {\it invariant} for $(\tau_t)$ if $\phi_0 \tau_t(x) = \phi_0(x)$ for all $x \in
\cla_0$ and $t \ge 0$. The support $p$ of a normal invariant state is a sub-harmonic
projection and $\phi^p_0$, the restriction of $\phi_0$ to $\cla^p_0$ is
a faithful normal invariant state for $(\tau^p_t)$. Thus asymptotic properties 
( ergodic, mixing ) of the dynamics $(\cla_0,\tau_t,\phi_0)$ is well determined by 
the asymptotic properties (ergodic, mixing respectively ) of the reduced dynamics 
$(\cla^p_0,\tau^p_t,\phi^p_0)$ provided $y=1$. For more details we refer to [Mo1]. 

\vsp
In case $\phi_0$ is faithful, normal and invariant for $(\tau_t)$, we recall [Mo1] that
$\clg=\{x \in \cla_0: \tilde{\tau}_t \tau_t(x) = x,\; t \ge 0 \}$ is von-Neumann sub-algebra of $\clf= \{ x  \in \cla_0:
\tau_t(x^*)\tau_t(x) = \tau_t(x^*x),\tau_t(x)\tau_t(x^*) = \tau_t(xx^*)\; \forall t \ge 0 \} $ and
the equality $\clg=\IC$ is a sufficient condition for $\phi_0$ to be strong mixing for $(\tau_t)$. Since the
backward process [AM] is related with the forward process via an anti-unitary operator we note that $\phi_0$ is
strongly mixing for $(\tau_t)$ if and only if same hold for $(\tilde{\tau}_t)$. We can also
check this fact by exploring faithfulness of $\phi_0$ and the adjoint relation [OP]. Thus $\IC \subseteq \tilde{\clg} \subseteq
\tilde{\clf}$ and equality $\IC=\tilde{\clg}$ is also a sufficient condition for strong mixing where $\tilde{\clf}$ and
$\tilde{\clg}$ are von-Neumann algebras associated with $(\tilde{\tau}_t)$. Thus we find two competing criteria for strong mixing.
However it is straight forward whether $\clf = \tilde{\clf}$ or $\clg=\tilde{\clg}$. Since given a dynamics it is difficult to describe
$(\tilde{\tau}_t)$ explicitly and thus this criterion $\clg=\IC$ is rather non-transparent. We prove in section 2 
that $\clg=\{ x \in \clf: \tau_t\sigma_s(x) = \sigma_s \tau_t(x),\; \forall t \ge 0.\; s \in \!R \}$ where
$\sigma=(\sigma_s:\;s \in \!R)$ is the Tomita's modular auto-morphism group [BR,OP] associated with $\phi_0$.
So $\clg$ is the maximal von-Neumann sub-algebra of $\cla_0$, where $(\tau_t)$ is an $*$-endomorphism [Ar], 
invariant by the modular auto-morphism group $(\sigma_s)$. Moreover $\sigma_s(\clg)=\clg$ for all $s \in \!R$ and 
$\tilde{\tau}_t(\clg)=\clg$ for all $t \ge 0$. Thus by a theorem of Takesaki [OP], there exists a norm one projection 
$\IE_{\clg}$ from $\cla_0$ onto $\clg$ which preserves $\phi_0$ i.e. $\phi_0 \IE = \phi_0$. Exploring
the fact that $\tilde{\tau}_t(\clg)=\clg$, we also conclude that the conditional expectation $\IE_{\clg}$ commutes with
$(\tau_t)$. This enables us to prove that $(\cla_0,\tau_t,\phi_0)$ is ergodic (strongly mixing) if and only if $(\clg,\tau_t,\phi_0)$
is ergodic (strongly mixing). Though $\tau_t(\clg) \subseteq \clg$ for all $t \ge 0$, equality may not hold in general.
However we have $$\bigcap_{t \ge 0} \tau_t(\clg)=\bigcap_{t \ge 0}\tilde{\tau}_t(\tilde{\clg})$$
where $\tilde{\clg}=\{ x \in \cla_0: \tau_t(\tilde{\tau}_t(x))=x,\; t \ge 0 \}$.
$\clg = \tilde{\clg}$ holds if and only if $\tau_t(\clg)=\clg,\;\tilde{\tau}_t(\tilde{\clg})=\tilde{\clg}$ for all $t \ge 0$. 
Thus $\clg_0= \bigcap_{t \ge 0} \tau_t(\clg)$ is the maximal von-Neumann sub-algebra invariant by the modular
automorphism so that $(\clg_0,\tau_t,\phi_0)$ is an $*-$automorphisms with $(\clg_0,\tilde{\tau}_t,\phi_0)$
as it's inverse dynamics. Once more there exists a conditional expectation $\IE_{\clg_0}:\cla_0 \raro \cla_0$
onto $\clg_0$ commuting with $(\tau_t)$. This ensures that $(\cla_0,\tau_t,\phi_0)$ is ergodic (strongly mixing)
if and only if $(\clg_0,\tau_t,\phi_0)$ is ergodic (strongly mixing). It is clear now that
$\clg_0=\tilde{\clg}_0$, thus $\clg_0=\IC$, a criterion for strong mixing, is symmetric or time-reversible. As an 
application in classical probability we can find an easy criteria for a stochastically complete Brownian flows [Mo5] on a 
Riemannian manifold driven by a family of complete vector fields to be strong mixing.   

\vsp
Exploring the criterion $\clg_0=\IC$ we also prove that for a type-I factor $\cla_0$ with center completely atomic,
strong mixing is equivalent to ergodicity when the time variable is continuous i.e. $\!R_+$ (Theorem 3.4). This result in particular
extends a result proved by Arveson [Ar] for type-I finite factor. In general, for discreet time dynamics $(\cla_0,\tau,\phi_0)$, 
ergodicity does not imply strong mixing property (not a surprise fact since we have many classical cases). We also prove that 
$\tau$ on a type-I von-Neumann algebra $\cla_0$ with completely atomic center is strong mixing 
if and only if it is ergodic and the point spectrum of $\tau$ in the unit circle i.e. $\{w \in S^1: \tau(x)=wx 
\;\;\mbox{for some non zero}\;\; x \in \cla_0 \}$ is trivial. The last result in a sense gives a direct proof of
a result obtained in section 7 of [BJKW] without being involved with Popescu dilation.      

\vsp
In section 3 we consider the unique up to isomorphism minimal forward weak Markov [AM,Mo1,Mo4] stationary process $\{ 
j_t(x),\; t \in \IT,\; x \in \cla_0 \}$ associated with $(\cla_0,\tau_t,\phi_0)$. We set a family of isomorphic von-Neumann 
algebras $\{ \cla_{[t}: t \in \IT \}$ generated by the forward process so that $\cla_{[t} \subseteq \cla_{[s}$ whenever 
$s \le t$. In this framework we construct a unique modulo unitary equivalence minimal dilation $(\cla_{[0},\alpha_t,\; t \ge 0,\phi)$, 
where $\alpha=(\alpha_t:t \ge 0)$ is a semigroup of $*-$endomorphism on a von-Neumann algebra $\cla_{[0}$ acting on a Hilbert space 
$\clh_{[0}$ with a normal invariant state $\phi$ and a projection $P$ in $\cla_{[0}$ so that

\NI (a) $P\cla_{[0}P=\pi(\cla_0)''$; 

\NI (b) $\Omega \in \clh_{[0}$ is a unit vector so that $\phi(X)=<\Omega, X \Omega>$;

\NI (b) $P\alpha_t(X)P=\pi(\tau_t(PXP))$ for $t \ge 0,\;X \in \cla_{[0}$; 

\NI (c) $ \{ \alpha_{t_n}(PX_nP) .....\alpha_{t_3}(PX_3P)\alpha_{t_2}(PX_2P)\alpha_{t_1}(PX_1P)\Omega:\; 0 \le t_1 \le t_2 ..\le t_n
,\; n \ge 1 \}, X_i \in \cla_{[0} \}$ is total in $\clh_{[0}$,

\NI where $\pi$ is the GNS representation of $\cla_0$ associated with the state $\phi_0$. In case $\phi_0$ is also faithful, we 
consider the backward process $(j^b_t)$ defined in [AM] associate with the KMS adjoint Markov semigroup and prove that commutant 
of $\cla_{[t}$ is equal to $\cla^b_{t]}= \{ j^b_s(x):\;x \in \cla_0, s \le t \}''$ for any fix $t \in \IT$. 

\vsp
As an application of our result on asymptotic behavior of a Markov semigroup, we also study a family of endomorphism 
$(\clb,\lambda_t)$ on a von-Neumann algebra. Following Powers [Po2] an endomorphism $\alpha_t: \clb_0 \raro \clb_0$ is called 
shift if $\bigcap_{t \ge 0}\alpha_t(\clb)$ is trivial. In general such a shift may not admit an invariant state [BJP]. Here 
we assume that $\lambda_t$ admits an invariant state $\psi$ and address how the shift property is related with Kolmogorov's 
property of the canonical Markov semigroup $(\cla_0,\tau_t,\psi_)$ on the support projection on the von-Neumann algebra 
$\pi_{\psi}(\clb)''$ of the state vector state in the GNS space $(\clh_{\pi},\pi,\Omega)$ associated with $(\clb,\psi)$. 
As a first step here we prove that Powers's shift property is equivalent to Kolmogorov's property of the adjoint Markov 
semigroup $(\tilde{\tau}_t)$. However in the last section we show that Kolmogorov's property of a Markov semigroup need not 
be equivalent to Kolmogorov's property of the KMS adjoint Markov semigroup. Thus Powers's shift property in general is not 
equivalent to Kolmogorov's property of the associated Markov semigroup.  

\vsp
Section 4 includes the main mathematical result by proving a criteria for the inductive limit state, associated with an 
invariant state of an injective endomorphism on a $C^*$ algebra, to be pure. To that end we explore the minimal weak 
Markov process associated with the reduced Markov semigroup on the corner algebra of the support projection and 
prove that the inductive limit state is pure if the Markov semigroup satisfies Kolmogorov's property. Further for
a lattice symmetric factor state, Kolmogorov's property is also necessary for purity of the inductive limit state.  

\vsp
The last section deals with an application of our main results on translation invariant state on quantum spin chain. We give a 
simple criteria for such a factor state to be pure and find its relation with Kolmogorov's property. Here we also deal with the 
unique temperature state i.e. KMS state on Cuntz algebra to illustrate that Powers shift property is not equivalent to Kolmogorov's
property of the associated canonical Markov map on the support projection. In fact this shows that Kolmogorov's property
is an appropriate notion to describe purity of the inductive state.

\newsection{ Time-reverse Markov semigroup and asymptotic properties: }

\vsp
In this section we will deal will a von-Neumann algebra $\cla$ and a completely
positive map $\tau$ or a semigroup $\tau=(\tau_t,\;t \ge 0 \}$ of such maps on $\cla$. 
We assume further that there exists a normal invariant state $\phi_0$ for $\tau$ and 
aim to investigate asymptotic properties of the Markov map. We say $(\cla_0,\tau_t,\phi_0)$ is
ergodic if $\{x: \tau_t(x)=x,\; t \ge 0 \}= \{zI, z \in \IC \}$ and we say mixing if 
$\tau_t(x) \raro \phi_0(x)$ in the weak$^*$ topology as $t \raro \infty$ for all $x \in \cla_0$. 

\vsp
For the time being we assume $\phi_0$ is faithful and recall following [OP,AM], the unique 
Markov map $\tilde{\tau}$ on $\cla_0$ which satisfies the following adjoint relation
\be
\phi_0(\sigma_{1/2}(x)\tau(y))=\phi_0(\tilde{\tau}(x)\sigma_{-1/2}(y))
\ee
for all $x,y \in \cla_0$ analytic elements for the Tomita's modular 
automorphism $(\sigma_t:\; t \in \IR)$ associated with a faithful normal 
invariant state for a Markov map $\tau$ on $\cla_0$. For more details 
we refer to the monograph [OP]. We also quote now [OP, Proposition 8.4 ] the 
following proposition without a proof. 

\vsp
\NI {\bf PROPOSITION 2.1: } Let $\tau$ be an unital completely positive normal
maps on a von-Neumann algebra $\cla_0$ and $\phi_0$ be a faithful normal 
invariant state for $\tau$. Then the following conditions are equivalent 
for $x \in \cla_0$:

\NI (a) $\tau(x^*x)=\tau(x^*)\tau(x)$ and $\sigma_s(\tau(x))= \tau(\sigma_s(x)),\; \forall \; s \in \!R;$

\NI (b) $\tilde{\tau} \tau (x)=x.$  

\NI Moreover $\tau$ restricted to the sub-algebra $\{x: \tilde{\tau}\tau(x) 
=x \}$ is an isomorphism onto the sub-algebra $\{x \in \cla_0: \tau\tilde{\tau}
(x) =x \}$ where $(\sigma_s)$ be the modular automorphism on $\cla_0$ 
associated with $\phi_0$. 

\vsp
In the following we investigate the situation further.

\vsp
\NI {\bf PROPOSITION 2.2: } Let $(\cla_0,\tau_t,\phi_0)$ be a quantum dynamical system and $\phi_0$ be
faithful invariant normal state for $(\tau_t)$. Then the following hold:

\NI (a) $\clg = \{x \in \cla_0: \tau_t(x^*x)=\tau_t(x^*) \tau_t(x),\; \tau_t(xx^*) = \tau_t(x)\tau_t(x^*),\; \sigma_s(\tau_t(x))
        = \tau_t (\sigma_s(x)),\; \forall \; s \in \!R,\; t \ge 0 \}$ and $\clg$ is $\sigma= (\sigma_s:\;s \in \!R)$ invariant
        and commuting with $\tau=(\tau_t:t \ge 0)$ on $\clg$. Moreover for all $t \ge 0,\; \tilde{\tau}_t(\clg)=\clg$ and the
        conditional expectation $\!E_{\clg}: \cla_0 \raro \cla_0$ onto $\clg_0$ commutes with $(\tau_t)$.

\NI (b) There exists a unique maximal von-Neumann algebra $\clg_0 \subseteq \clg \bigcap \tilde{\clg}$ so that 
$\sigma_t(\clg_0)=\clg_0$ for all $t \in \!R$ and $(\clg_0,\tau_t,\phi_0)$ is an automorphism where for any $t \ge 0$, 
$\tilde{\tau}_t\tau_t=\tau_t\tilde{\tau}_t=1$ on $\clg_0$. Moreover the conditional expectation $\!E_{\clg_0}:\cla_0 \raro \cla_0$ 
onto $\clg_0$ commutes with $(\tau_t)$ and $(\tilde{\tau}_t)$.

\vsp
\NI {\bf PROOF:} The first part of (a) is a trivial consequence of Proposition 2.1
once we note that $\clg$ is closed under the action $x \raro x^*$. For the
second part we recall [Mo1] that $\phi_0(x^*JxJ)-
\phi_0(\tau_t(x^*)J\tau_t(x)J)$ is monotonically increasing with $t$ and
thus for any fix $t \ge 0$ if $\tilde{\tau}_t\tau_t(x)=x$ then
$\tilde{\tau}_s\tau_s(x)=x$ for all $0 \le s \le t$. So the
sequence $\clg_t=\{ x \in \cla_0:\tilde{\tau}_t\tau_t(x)=x \}$ of von-Neumann
sub-algebras decreases to $\clg$ as $t$ increases to $\infty$ i.e.
$\clg = \bigcap_{t \ge 0} \clg_t$. Similarly we also have $\tilde{\clg}=
\bigcap_{ t \ge 0} \tilde{\clg}_t$.

\vsp
Since $\tilde{\clg}_t$ monotonically decreases to $\tilde{\clg}$ as $t $ increases to infinity
for any $s \ge 0$ we claim that $\tau_s(\tilde{\clg}^1) = \bigcap_{t \ge 0} \tau_s(\tilde{\clg}_t^1)$, where we 
have used the symbol $\cla^1=\{x \in \cla: ||x||=1 \}$ for a von-Neumann algebra $\cla$. We will prove 
the non-trivial inclusion. To that end let $x \in \bigcap_{t \ge 0} \tau_s(\tilde{\clg}_t^1)$ i.e. for each $t \ge 0$ 
there exists $y_t \in \tilde{\clg}_t^1$ so that $\tau_s(y_t)=x$. By weak$^*$ compactness of the unit ball of $\cla_0$,
we extract a subsequence $t_n \raro \infty$ so that $y_{t_n} \raro y$ as $t_n \raro \infty$ for some $y \in \cla_0$. 
The von-Neumann algebras $\tilde{\clg}_t$ being monotonically decreasing, for each $m \ge 1$, $y_{t_n} \in \tilde{\clg}_{t_m}$
for all $n \ge m$. $\tilde{\clg}_{t_m}$ being a von-Neumann algebra, we get $y \in \tilde{\clg}_{t_m}$. As this holds for each 
$m \ge 1$, we get $y \in \tilde{\clg}$. However by normality of the map $\tau_s$, we also have $x=\tau_s(y)$. Hence $x \in 
\tau_s(\tilde{\clg}^1)$.    

\vsp
Now we verify that $\bigcap_{s \ge r } \tau_s(\tilde{\clg}^1) = \bigcap_{ s \ge r } \bigcap_{t \ge 0} \tau_{s+t}(\clg^1_t)=
\bigcap_{t \ge 0} \bigcap_{s \ge r} \tau_{s+t}(\clg^1_t) = \bigcap_{t \ge r} \bigcap_{0 \le s \le t}\tau_t(\clg^1_s)$,
where we have used $\tau_t(\clg^1_t)=\tilde{\clg}^1_t$ being isomorphic. Since $\clg_t$
are monotonically decreasing with $t$ we also note that
$\bigcap_{0 \le s \le t} \tau_t(\clg^1_s) = \tau_t(\clg^1_t)$. Hence for any $r \ge 0$
\be
\bigcap_{s \ge r } \tau_s(\tilde{\clg}^1) = \tilde{\clg}^1
\ee
From (2.2) with $r=0$ we get $\tilde{\clg}^1 \subseteq \tau_t(\tilde{\clg}^1)$ for all $t \ge 0$. For any $t \ge 0$
we also have $\tau_t(\tilde{\clg}^1) \subseteq \bigcap_{s \ge t } \tau_s(\tilde{\clg}^1) = \tilde{\clg}^1$. Hence we
conclude $\tau_t(\tilde{\clg}^1) = \tilde{\clg}^1$ for any $t \ge 0$. Now we can easily remove the restriction to show 
that $\tau_t(\tilde{\clg}) = \tilde{\clg}$ for any $t \ge 0$ by linearity. By symmetry $\tilde{\tau}_t(\clg)=\clg$ for 
any $t \ge 0$.

\vsp
Since $\clg$ is invariant under the modular automorphism $(\sigma_s)$ by a theorem of
Takesaki [AC] there exists a norm one projection $\!E_{\clg}:\cla \raro \cla$ with range equal
to $\clg$. We claim that $\!E_{\clg}$ commutes with $(\tau_t)$. To that end we verify for any $x \in \cla_0$ and
$y \in \clg$ the following equalities:

$$ <J_{\clg}yJ_{\clg}\omega_0, \!E_{\clg}(\tau_t(x)) \omega_0> = <J_0yJ_0\omega_0, \tau_t(x) \omega_0>$$
$$=<J_0\tilde{\tau}_t(y)J_0 \omega_0,x \omega_0> = <J_{\clg}\tilde{\tau}_t(y)J_{\clg} \omega_0, \!E_{\clg}(x)
\omega_0> $$
$$ = <J_{\clg}yJ_{\clg}\omega_0,\tau_t(\!E_{\clg}(x))\omega_0)> $$
where we used the fact that $\tilde{\tau}(\clg)=\clg$ for the third equality and range of $\IE_{\clg}$ is indeed $\clg$ is
used for the last equality. This completes the proof of (a).

\vsp
Now for any $s \ge 0$, it is obvious that $\tilde{\tau}_s(\tilde{\clg}) \subseteq 
\bigcap_{t \ge s} \tilde{\tau}_s(\tilde{\clg}_t)$. In the following we prove equality in the above relation. Let 
$x \in \bigcap_{t \ge s} \tilde{\tau}_s(\tilde{\clg}_t)$ i.e. there exists elements $y_t \in \tilde{\clg}_t$ so that 
$x=\tilde{\tau}_s(y_t)$ for all $t \ge s$. If so then we have $\tau_s(x)=y_t$ for all $t \ge s$ as $\tilde{\clg}_t \subseteq 
\tilde{\clg}_s$. Thus for any $t \ge s$, $y_t=y_s \in \tilde{\clg}$ and $x \in \tilde{\tau}_s(\tilde{\clg})$.    

\vsp
Now we verify the following elementary relations:
$\tilde{\tau}_s(\tilde{\clg}) = \bigcap_{t \ge s} \tilde{\tau}_s
\tau_t(\clg_t)= \bigcap_{t \ge s}\tilde{\tau}_s \tau_s(\tau_{t-s}(\clg_t)))= \bigcap_{t \ge s} \tau_{t-s}(\clg_t)=
\bigcap_{t \ge 0}\tau_t(\clg_{s+t})$ where we have used the fact that $\tau_{t-s}(\clg_t) \subseteq \clg_s$.
Thus we have $\bigcap_{s \ge 0}\tau_s(\clg) \subseteq  \bigcap_{s \ge 0}\tilde{\tau}_s(\tilde{\clg})$. By the dual
symmetry, we conclude the reverse inclusion and hence
\be
\bigcap_{s \ge 0}\tau_s(\clg) =  \bigcap_{s \ge 0}\tilde{\tau}_s(\tilde{\clg})
\ee

\vsp
We set von-Neumann algebra $\clg_0 = \bigcap_{s \ge 0}\tau_s(\clg)$. Thus $\clg_0 \subseteq \clg$ and also $\clg_0 \subseteq 
\tilde{\clg}$ by (2.3) and for each $t \ge 0$ we have $\tau_t\tilde{\tau}_t=\tilde{\tau}_t \tau_t=1$ on $\clg_0$. Since 
$\tau_s(\clg)$ is monotonically decreasing, we also note that $\tau_t(\clg_0)= \bigcap_{s \ge 0}\tau_{s+t}(\clg)=\clg_0$. 
Similarly $\tilde{\tau}_t(\clg_0)=\clg_0$ by (2.3). That $\clg_0$ is invariant by the modular group $\sigma$ follows since 
$\clg$ is invariant by $\sigma=(\sigma_t)$ which is commuting with $\tau=(\tau_t)$ on $\clg$. Same is also true for 
$(\tilde{\tau}_t)$ by (2.3). By Takesaki's theorem [AC] once more we guarantee that there exists a conditional expectation 
$\!E_{\clg_0}: \cla_0 \raro \cla_0$ with range equal to $\clg_0$. Since $\tilde{\tau_t}(\clg_0)= \clg_0$, once more by repeating 
the above argument we conclude that $E_{\clg_0} \tau_t= \tau_t E_{\clg_0}$ on $\cla_0$. By symmetry of the argument, $\!E_{\clg_0}$ 
is also commuting with $\tilde{\tau}=(\tilde{\tau}_t)$ \qed

\vsp
We have the following reduction theorem.  

\vsp
\NI {\bf THEOREM 2.3: } Let $(\cla_0,\tau_t,\phi_0)$ be as in Proposition 2.2. 
Then the following statements are equivalent:

\NI (a) $(\cla_0,\tau_t,\phi_0)$ is mixing ( ergodic );

\NI (b) $(\clg, \tau_t,\phi_0)$ is mixing ( ergodic );

\NI (c) $(\clg_0,\tau_t,\phi_0)$ is mixing ( ergodic ). 

\vsp
\NI {\bf PROOF: } That (a) implies (b) is obvious. By Proposition 2.2. we have $\!E_{\clg}\tau_t(x) = \tau_t \!E_{\clg}(x)$ for any 
$x \in \cla_0$ and $t \ge 0$. Fix any $x \in \cla_0$. Let $x_{\infty}$ be any weak$^*$ limit point of the net $\tau_t(x)$ as $t 
\raro \infty$ which is an element in $\clg$ [Mo1]. In case (b) is true, we find that $x_{\infty} = \!E_{\clg}(x_{\infty}) = 
\phi_0(\!E_{\clg}(x))=\phi_0(x)1$. Thus $\phi_0(x)1$ is the unique limit point, hence weak$^*$ limit of $\tau_t(x)$ as 
$t \raro \infty$ is $\phi_0(x)1$. The equivalence statement for ergodicity also follows along the same line since the conditional 
expectation $\!E_{\cli}$ on the the von-Neumann algebra $\cli = \{ x: \tau_t(x)=x,\; t \ge 0 \}$ commutes with $(\tau_t)$ and thus 
satisfies $\!E_{\cli} \!E_{\clg} = \!E_{\clg} \!E_{\cli} = \!E_{\cli}$. This completes the proof that (a) and (b) are equivalent. 
That (b) and (c) are equivalent follows essentially along the same line since once more there exists a conditional expectation 
from $\clg$ to $\clg_0$ commuting with $(\tau_t)$ and any weak$^*$ limit point of the net $\tau_t(x)$ as $t$ diverges to infinity 
belongs to $\tau_s(\clg)$ for each $s \ge 0$, thus in $\clg_0$. We omit the details. \qed

\vsp
Now we investigate asymptotic behavior for quantum dynamical system dropping the assumption that $\phi_0$ is faithful. Let
$p$ be the support projection of the normal state $\phi_0$ in $\cla_0$. Thus we have $\phi_0(p\tau_t(1-p)p)=0$ for 
all $t \ge 0$, $p$ being the support projection we have $p\tau_t(1-p)p=0$ i.e. $p$ is a sub-harmonic projection 
in $\cla_0$ for $(\tau_t)$ i.e.  $\tau_t(p) \ge p$ for all $t \ge 0$. Then it is simple to check that 
$(\cla_0^p,\tau^p_t,\phi^p_0)$ is a quantum dynamical semigroup where $\cla_0^p=p\cla_0p$ and 
$\tau^p_t(x)=p\tau_t(pxp)p$ for $x \in \cla^p_0$ and $\phi^p_0(x)=\phi_0(pxp)$ is faithful on $\cla^p_0$. 
In Theorem 3.6 and Theorem 3.12 in [Mo1] we have explored how ergodicity and strong mixing of the original dynamics 
$(\cla_0,\tau_t,\phi_0)$ can be determined by that of the reduced dynamics $(\cla_0^p,\tau^p_t,\phi^p_0)$. Here we add 
one more result in that line of investigation.

\vsp
\NI {\bf THEOREM 2.4:} Let $(\cla_0,\tau_t,\phi_0)$ be a quantum dynamical systems with a normal invariant state $\phi_0$
and $p$ be a sub-harmonic projection for $(\tau_t)$. If $\mbox{s-limit}_{t \raro \infty}\tau_t(p)=1$
then the following statements are equivalent:

\NI (a) $||\phi\tau_t-\phi_0|| \raro 0$ as $t \raro \infty$ for any normal state on $\phi$ on $\cla_0$.

\NI (b) $||\phi^p \tau^p_t- \phi^p_0|| \raro 0$ as $t \raro \infty$ for any normal state $\phi^p$ on $\cla^p_0$.

\vsp
\NI {\bf PROOF: } That (a) implies (b) is trivial. For the converse we write
$||\phi \tau_t-\phi_0|| = \mbox{sup}_{x:||x|| \le 1} |\phi \tau_t(x)-\phi_0(x)|
\le \mbox{sup}_{\{x:||x|| \le 1\}} |\phi \tau_t(pxp)-\phi_0(pxp)| +
\mbox{sup}_{ \{x:||x|| \le 1 \} } |\phi \tau_t(pxp^{\perp})| +
\mbox{sup}_{ \{x:||x|| \le 1 \} } |\phi \tau_t(p^{\perp}xp)| +
\mbox{sup}_{ \{x:||x|| \le 1 \} } |\phi \tau_t(p^{\perp}xp^{\perp})|$.
Since $\tau_t((1-p)x) \raro 0$ in the weak$^*$ topology and
$|\phi \tau_t(xp^{\perp})|^2 \le |\phi \tau_t(xx^*)|\phi(\tau_t(p^{\perp}))| \le ||x||^2 \phi(\tau_t(p^{\perp})$
it is good enough if we verify that (a) is
equivalent to  $\mbox{sup}_{\{ x:||x|| \le 1 \} } |\phi\tau_t(pxp)-\phi_0(pxp)| \raro 0$
as $t \raro \infty$. To that end we first note that limsup$_{t \raro
\infty} \mbox{sup}_{x:||x|| \le 1} |\psi(\tau_{s+t}(pxp))-\phi_0(pxp)|$ is independent
of $s \ge 0$ we choose. On the
other hand we write $\tau_{s+t}(pxp) = \tau_s(p\tau_t(pxp)p) +
\tau_s(p\tau_t(pxp)p^{\perp})+ \tau_s(p^{\perp}\tau_t(pxp)p) +
\tau_s(p^{\perp}\tau_t(pxp)p^{\perp})$ and use the fact for any
normal state $\phi$ we have
$\mbox{limsup}_{t \raro \infty}\mbox{sup}_{x:||x|| \le 1}
|\psi(\tau_s(z\tau_t(pxp)p^{\perp})|\le ||z||\;|\psi(\tau_s(p^{\perp}))| $ for all $z \in \cla_0$.
Thus by our hypothesis on the support projection we conclude that (a) hold whenever (b) is true. \qed

\vsp
In case the time variable is continuous and the von-Neumann algebra is the set of bounded linear operators on a finite dimensional 
Hilbert space $\clh_0$, by exploring Lindblad's representation [Li], Arveson [Ar] shows that a quantum dynamical semigroup with a 
faithful normal invariant state is ergodic if and only if the dynamics is mixing. In the following we prove a more general result 
exploring the criteria that we have obtained in Theorem 2.3. Note at this point that we don't even need the generator of the Markov 
semigroup to be a bounded operator for which Lindblad's representation is not yet understood with full generality [CE].  

\vsp
\NI {\bf THEOREM 2.5:} Let $\cla_0$ be type-I with center completely atomic and $(\tau_t:t \in \!R)$ 
admits a normal invariant state $\phi_0$. Then $(\cla_0,\tau_t,\phi_0)$ is strong mixing if and only if 
$(\cla_0,\tau_t,\phi_0)$ is ergodic. 

\vsp
\NI {\bf PROOF: } We first assume that $\phi_0$ is also faithful. We will verify now the criteria that 
$\clg_0$ is trivial when $(\tau_t)$ is ergodic. Since $\clg_0$ is invariant by the modular 
auto-morphism group associated with the faithful normal state $\phi_0$, by a theorem of Takesaki [Ta] there 
exists a faithful normal norm one projection from $\cla_0$ onto $\clg_0$. Now since $\cla_0$ is a von-Neumann 
algebra of type-I with center completely atomic, a result of E. Stormer [So] says that 
$\clg_0$ is also type-I with center completely atomic. 

\vsp
Let $Q$ be a central projection in $\clg_0$. Since $\tau_t(Q)$ is also a central projection and $\tau_t(Q) \raro Q$ as 
$t \raro 0$ we conclude that $\tau_t(Q)=Q$ for all $t \ge 0$ (center of $\clg$ being completely atomic and time 
variable $t$ is continuous ). Hence by ergodicity we conclude that $Q = 0$ or $1$. Hence $\clg_0$ can be identified 
with $\clb(\clk)$ for a separable Hilbert space $\clk$. Since $(\tau_t)$ on $\clb(\clk)$ is an automorphism we find a 
self-adjoint operator $H$ in $\clk$ so that $\tau_t(x)=e^{itH}xe^{-itH}$ for any $x \in \clb(\clk)$. Since it admits 
an ergodic faithful normal state, by [Fr, Mo1] we conclude that $\{ x \in \clb(\clk): x e^{itH} = e^{itH} x,\; t \in \!R \} 
= \IC$, which holds if and only if $\clk$ is one dimensional. Hence $\clg_0 = \IC$. 

\vsp
Now we deal with the general situation. Let $p$ be the support projection of $\phi_0$ in $\cla_0$ and $\cla_0$ being 
a type-I von-Neumann algebra with centre completely atomic, the center of $\cla_0^p=p\cla_0p$ being equal to the corner 
of the center of $\cla_0$ i.e. $p\cla_0 \bigcap \cla'_0p$, is  also a type-I von-Neumann algebra with completely atomic 
centre. $(\cla_0,\tau_t,\phi_0)$ being ergodic, we have $\tau_t(p) \uparrow 1$ as $t \uparrow \infty$ in the weak$^*$ topology 
and $(\cla^p_0,\tau_t^p,\phi^p_0)$ is ergodic. Thus by the first part of the argument, $(\cla_0^p,\tau_t^p,\phi^p_0)$ is strongly 
mixing. Hence by Theorem 3.12 in [Mo1] we conclude that $(\cla_0,\tau_t,\phi_0)$ is also strong mixing. This completes the proof. \qed

\vsp
We end this section with another simple application of Theorem 2.3 by proving a result originated in [FNW1,FNW2,BJKW]. 

\vsp
\NI {\bf THEOREM 2.6: } Let $\cla_0$ be a type-I von-Neumann algebra with center completely atomic and $\tau$ be 
a completely positive map with a faithful normal invariant state $\phi_0$. Then the following are equivalent:

\NI (a) $(\cla_0,\tau_n,\phi_0)$ is strong mixing.

\NI (b) $(\cla_0,\tau_n,\phi_0)$ is ergodic and $\{ w \in S^1,\;\tau(x)=wx, \mbox{ for some non zero} x \in \cla_0 \} = 
\{ 1 \}$, where $S^1=\{w \in \IC: |w|=1 \}.$    

\vsp
\NI {\bf PROOF: } That `(a) implies (b)' is rather simple and true in general for any von-Neumann algebra. To that end let 
$\tau(x)=wx$ for some $x \ne 0$ and $|w|=1$. Then $\tau^n(x)=w^nx$ and since the sequence $w^n$ has a limit point say $z, |z|=1$ 
we conclude by strong mixing that $zx=\phi_0(x)I$. Hence $x$ is a scaler and thus $x=\tau(x),\;x \ne 0$. So $w=1$ and $x=\phi_0(x)I$. 
Ergodic property also follows by strong mixing as $x=\phi_0(x)I$ for any $x$ for which $\tau(x)=x$. 

\vsp
Now for the converse we will use our hypothesis that $\phi_0$ is faithful and $\cla_0$ is a type-I von-Neumann algebra with
completely atomic. To that end we plan to verify that $\clg_0$ consists of scalers only and appeal to Theorem 2.3 for strong mixing. 
Since there exists a conditional expectation from $\cla_0$ onto $\clg_0$, by a Theorem of Stormer [So] $\clg_0$ is once more a type-I 
von-Neumann algebra with center completely atomic. Let $E$ be a non-zero atomic projection in the center of $\clg_0$. $\tau$ being an 
automorphism on $\clg_0$, each element in the sequence $\{ \tau_k(E): k \ge 0 \}$ is an atomic projection in the 
center of $\clg_0$. If $\tau_n(E) \bigcap \tau_m(E) \ne 0$ and $n \ge m$ we find that $ \tau_m(\tau_{n-m}(E) \bigcap E) \neq 0$ 
and thus by faithful and invariance property of $\phi_0$, we get $\phi(\tau_{n-m}(E) \bigcap E) > 0$. Once more by 
faithfulness we find that $\tau_{n-m}(E) \bigcap E \ne 0$. So by atomic property of $E$ and $\tau_{n-m}(E)$ 
we conclude that $\tau_{n-m}(E)=E$. Thus either the elements in the infinite sequence $E, \tau(E),...., \tau^n(E)....$ 
are all mutually orthogonal or there exists a least positive integer $n \ge 1$ so that the projections 
$E,\tau(E),.., \tau_{n-1}(E)$ are mutually orthogonal and $\tau^n(E)=E$. However for such an infinite sequence with mutually 
orthogonal projection we have $1 = \phi_0(I) \le \phi_0(\bigcup_{0 \le n \le m-1} \tau_n(E) ) = m \phi_0(E)$ for all $m \ge 1$. 
Hence $\phi_0(E)=0$ contradicting that $E$ is non-zero and $\phi_0$ is faithful.  

\vsp
Thus for any $w \in S^1$ with $w^n=1$, we have $\tau(x)=wx,$ where $x=\sum_{0 \le k \le n-1}w^k\tau_k(E) \ne 0 $. 
By (b) we have $w=1$. Hence $n=1$. In other words we have $\tau(E)=E$ for any atomic projection in the center of
$\clg_0$. Now by ergodicity we have $E=I$. Thus $\clg_0$ is a type-I factor say isomorphic to $\clb(\clk)$ for 
some Hilbert space $\clk$ and $\tau(x)=uxu^*$ for some unitary element $u$ in $\clg_0$. Since $(\clg_0,\tau_n,\phi_0)$ 
is ergodic we have $\{u,u* \}''= \clb(\clk)$, which holds if and only if $\clk$ is one dimensional 
( check for an alternative proof that $\tau(u)=u$, thus $u=I$ by ergodicity and thus $\tau(x)=x$ for all $x \in \clg_0$ ). 
Hence $\clg_0=\IC$. This complete the proof that (b) implies (a). \qed

\newsection{ Minimal endomorphisms and Markov semigroups : }

\vsp
An E$_0$-semigroup $(\alpha_t)$ is a weak$^*$-continuous one-parameter semigroup
of unital $^*$-endomorphisms on a von-Neumann algebra $\cla$ acting on a Hilbert
space $\clh$. Following [Po1,Po2,Ar] we say $(\alpha_t)$ is a {\it shift } if $\bigcap_{t \ge 0}
\alpha_t(\cla)=\IC$. For each $t \ge 0$, $\alpha_t$ being an endomorphism, $\alpha_t(\cla)$ is 
itself a von-Neumann algebra and thus $\bigcap_{t \ge 0}\alpha_t(\cla)$ is a limit of 
a sequence of decreasing von-Neumann algebras. Exploring this property Arveson proved that 
$(\alpha_t)$ is pure if and only if $||\psi_1\alpha_t-\psi_2\alpha_t|| \raro 0$ as $t \raro 
\infty$ for any two normal states $\psi_1,\psi_2$ on $\cla$. These criteria gets further simplified 
in case $(\alpha_t)$ admits a normal invariant state $\psi_0$ for which we have $(\alpha_t)$ 
is a shift (in his terminology it is called pure, here we prefer Powers's terminology as the last section will illustrate 
a shift need not be pure in its inductive limit ) if and only if $||\psi\alpha_t-\psi_0|| \raro 0$ as $ t \raro \infty$ 
for any normal state $\psi$. In such a case $\psi_0$ is the unique normal invariant state.  However a shift $(\alpha_t)$ 
in general may not admit a normal invariant state [Po2,BJP] and this issue is itself an interesting problem.  

\vsp 
One natural question that we wish to address here whether similar result is also true for a Markov 
semigroup $(\tau_t)$ defined on an arbitrary von-Neumann algebra $\cla_0$. This issue is already 
investigated in [Ar] where $\cla_0=\clb(\clh)$ and the semigroup $(\tau_t)$ is assumed to be 
continuous in the strong operator topology. He explored associated minimal dilation to an E$_0$-semigroups 
and thus make possible to prove that the associated $E_0$-semigroup is a shift if and only if 
$||\phi_1\tau_t-\phi_2\tau_t|| \raro 0$ as $t  \raro \infty$ for any two normal states 
$\phi_1,\phi_2$ on $\cla_0$. In case $(\tau_t)$ admits a normal invariant state the 
criteria gets simplified once more. In this section we will investigate this issue further for an arbitrary 
von-Neumann algebra assuming that $(\tau_t)$ is admits a normal invariant state $\phi_0$. 

\vsp
To that end, we consider [Mo1] the minimal stationary weak Markov forward process 
$(\clh,F_{t]},j_t,\Omega,\; t \in \!R)$ and Markov shift $(S_t)$ associated with 
$(\cla_0,\tau_t,\phi_0)$ and set $\cla_{[t}$ to be the von-Neumann algebra generated 
by the family of operators $\{j_s(x): t \le s < \infty,\; x \in \cla_0 \}$. We recall 
that $j_{s+t}(x)=S^*_tj_s(x)S_t,\; t,s \in \!R$ and thus $\alpha_t(\cla_{[0}) \subseteq \cla_{[0}$ 
whenever $t \ge 0$. Hence $(\alpha_t,\; t \ge 0)$ is a E$_0$-semigroup on $\cla_{[0}$ with a 
invariant normal state $\Omega$ and 
$$ j_s(\tau_{t-s}(x))=F_{s]}\alpha_t(j_{t-s}(x))F_{s]} \eqno (3.1) $$ 
for all $x \in \cla_0$. We consider the GNS Hilbert space $(\clh_{ \pi_{\phi_0} }, \pi_{\phi_0}(\cla_0),\omega_0)$ 
associated with $(\cla_0,\phi_0)$ and define a Markov semigroup $(\tau_t^{\pi})$ on 
$\pi(\cla_0)$ by $\tau^{\pi}_t(\pi(x))= \pi(\tau_t(x)$. Furthermore we now identify $\clh_{\phi_0}$ as the 
subspace of $\clh$ by the prescription $\pi_{\phi_0}(x)\omega_0 \raro j_0(x)\Omega$. In such a case $\pi(x)$ 
is identified as $j_0(x)$ and aim to verify for any $ t \ge 0$ that
$$ \tau^{\pi}_t(PXP)=P\alpha_t(X)P \eqno (3.2) $$ for all $X \in \cla_{[0}$
where $P$ is the projection from $\clh$ on the GNS space. We use induction on $n \ge 1$. If $X=j_s(x)$ for some 
$s \ge 0$, (4.2) follows from (4.1). Now we assume that (3.2) is true for any element of the form $j_{s_1}(x_1)...j_{s_n}
(x_n)$ for any $s_1,s_2,...,s_n \ge 0$ and $x_i \in \cla_0$ for $1 \le i \le n$. Fix any $s_1,s_2,,s_n,s_{n+1} 
\ge 0$ and consider $X=j_{s_1}(x_1)...j_{s_{n+1}}(x_{n+1})$. Thus 
$P\alpha_t(X)P=j_0(1)j_{s_1+t}(x_1)...j_{s_{n+t}}(x_{n+1})j_0(1)$. If $s_{n+1} \ge s_n$
we use (3.1) to conclude (3.2) by our induction hypothesis. Now suppose $ s_{n+1} \le s_n$. In that case if $s_{n-1} 
\le s_{n}$ we appeal to (3.1) and induction hypothesis to verify (3.2) for $X$. Thus we are left to consider 
the case where $s_{n+1} \le s_n \le s_{n-1}$ and by repeating this argument we are left to check only the case
where $s_{n+1} \le s_n \le s_{n-1} \le .. \le s_1 $. But $s_1 \ge 0=s_0$ thus we can appeal to (3.1) at the end of the
string and conclude that our claim is true for all elements in the $*-$ algebra generated by these elements 
of all order. Thus the result follows by von-Neumann density theorem.  We also 
note that $P=\tau^{\pi}_t(1)$ is a sub-harmonic projection [Mo1] for $(\alpha_t:t \ge 0)$ i.e. 
$\alpha_t(P) \ge P$ for all $t \ge 0$.  

\vsp
\NI {\bf PROPOSITION 3.1: } Let $(\cla_0,\tau_t,\phi_0)$ be a quantum dynamical semigroup with a normal invariant 
state for $(\tau_t)$. Then the GNS space $\clh_{\pi_{\phi_0}}$ associated with the normal state $\phi_0$ on $\cla_0 $ 
can be realized as a closed subspace of a unique Hilbert space $\clh_{[0}$ up to isomorphism so that the following hold:

\NI (a) There exists a von-Neumann algebra $\cla_{[0}$ acting on $\clh_{[0}$ and a unital $*$-endomorphism 
$(\alpha_t,\; t \ge 0 )$ on $\cla_{[0}$ with a pure vector state $\phi(X)=<\Omega,X \Omega>$, $\Omega \in \clh_{[0}$
invariant for $(\alpha_t:t \ge 0)$. 

\NI (b) $P \cla P$ is isomorphic with $\pi(\cla_0)$ where $P$ is the projection onto $\clh_{\pi_{\phi_0}}$; 

\NI (c) $P\alpha_t(X)P=\tau^{\pi}_t(PXP)$ for all $t \ge 0$ and $X \in \cla_{[0}$; 

\NI (d) The closed span generated by the vectors $\{ \alpha_{t_n}(PX_nP)....\alpha_{t_1}(PX_1P)\Omega: 0 \le t_1 \le t_2 \le ..
\le t_k \le ....t_n, X_1,..,X_n \in \cla_{[0}, n \ge 1 \}$ is $\clh_{[0}$. 

\vsp
\NI {\bf PROOF:} The uniqueness up to isomorphism follows from the minimality property (d). \qed 

\vsp
Following the literature [Vi,Sa,BhP] on dilation we say $(\cla_{[0},\alpha_t,\phi)$ is the minimal E$_0$-semigroup associated 
with $(\cla_0,\tau_t,\phi_0)$. By a theorem [Ar, Proposition 1.1 ] we conclude that $\bigcap_{t \ge 0}\alpha_t(\cla_{[0}) = 
\IC$ if and only if for any normal state $\psi$ on $\cla_{[0}$, $||\psi\alpha_t-\psi_0|| \raro 0$ as $t \raro  \infty$, where 
$\psi_0(X)=<\Omega,X\Omega>$ for $X \in \cla_{0]}$. In the following proposition we explore that fact that $P$ is a 
sub-harmonic projection for $(\alpha_t)$ and by our construction $\alpha_t(P)=F_{t]} \uparrow 1$ as $t \raro \infty$. 

\vsp
\NI {\bf PROPOSITION 3.2:} $||\psi\alpha_t-\psi_0|| \raro 0$ as $t \raro \infty$ for all normal state
$\psi$ on $\cla_{[0}$  if and only if $||\phi\tau_t-\phi_0|| \raro 0$ as $t \raro \infty$ for all normal state
$\phi$ on $\pi(\cla_0)''$, where $\pi$ is the GNS space associated with $(\cla_0,\phi_0)$. 
 
\vsp
\NI {\bf PROOF:} Since $F_{s]} \uparrow 1$ in strong operator topology by our construction and $\pi(\cla_0)$ is isomorphic to 
$F_{0]}\cla_{[0}F_{0]}$, we get the result by a simple application of Theorem 2.4. \qed   

\vsp
\NI {\bf THEOREM 3.3:} Let $\tau=(\tau_t,\;t \ge 0)$ be a weak$^*$ continuous Markov semigroup on $\cla_0$ 
with an invariant normal state $\phi_0$. Then there exists a weak$^*$ continuous E$_0$-semigroup 
$\alpha=(\alpha_t,\; t \ge 0)$ on a von-Neumann algebra $\cla_{[0}$ acting on a Hilbert space $\clh$ so that 
$$P\alpha_t(X)P=\tau^{\pi}_t(PXP),\; t \ge 0 $$ 
for all $X \in \cla_{[0}$, where $P$ is a sub-harmonic projection for $(\alpha_t)$ 
such that $\alpha_t(P) \uparrow I$.

\NI Moreover the following statements are equivalent:

\NI (a) $\bigcap_{t \ge 0}\alpha_t(\cla_{[0})=\!C$ 

\NI (b) $||\phi\tau_t^{\pi}-\phi_0|| \raro 0$ as $t \raro \infty$ for any normal state $\phi$ 
on $\pi(\cla_0)''$.

\vsp
\NI {\bf PROOF:} For convenience of notation we denote $\pi(\cla_0)''$ as $\cla_0$ in the following proof.
That (a) and (b) are equivalent follows by a Theorem of Arveson [Ar] and Proposition 3.2. \qed

\vsp
Following [AM,Mo1] we say $(\clh,S_t,F_{t]},\Omega)$ is a {\it Kolmogorov's shift } if strong 
$\mbox{lim}_{t \raro - \infty}F_{t]}=|\Omega><\Omega|$. We also recall here that Kolmogorov's shift property holds if and only if 
$\phi_0(\tau_t(x)\tau_t(y)) \raro  \phi_0(x)\phi_0(y)$ as $t \raro \infty$ for all $x,y \in \cla_0$. In such a case 
$\cla=\clb(\clh)$ ( see the paragraph before Theorem 3.9 in [Mo1] ). If $\phi_0$ is faithful then $\cla_0$ and $\pi(\cla_0)$ are 
isomorphic, thus $\bigcap_{t \ge 0}\alpha_t(\cla_{[0})=\!C$ if and only if $||\phi\tau_t-\phi_0|| \raro 0$ as $t \raro \infty$ 
for any normal state $\phi$ on $\cla_0$. Such a property is often called {\it strong ergodic property}. Our next result says that 
there is a duality between strong ergodicity and Kolmogorov's shift property. To that end  we recall the backward process 
$(\clh,j^b_t,F_{[t},\Omega)$ as defined in [AcM,Mo1] where $F_{t]}$ be the projection on the subspace generated by the vectors 
$\{\lambda: \IR \raro \cla_0: \mbox{ support of } \lambda \subseteq (-\infty, t] \}$ and for any $x \in \cla_0$, 
$j^b_t(x)$ is the trivial extension of it's action on $F_{t]}$ which takes an typical vector $\lambda$ to 
$\lambda'$ where $\lambda'(s)=\lambda(s)$ for any $s < t$ and $\lambda'(t)=\lambda(t)\sigma_{i \over 2}(x).$ For any analytic 
element $x$ for the automorphism group, we check first that $j^b_t$ is indeed an isometry if $x$ is so. Now we extend as 
analytic elements are weak$^*$ dense to all isometrics and extend by linearity to all elements of $\cla_0$. We recall here that 
we have backward Markov property for the process $(j^b_s)$ as $F_{[t}j^b_s(x)F_{[t}=j^b_t(\tilde{\tau}_{t-s}(x))$ for all
$t \ge s$ where $(\cla_0,\tilde{\tau}_t,\;t \ge 0,\phi_0$ is the dual Markov semigroup defined in (3.1). As in the forward process 
we have now $F_{[t}\cla^b_{t]}F_{[t}=j^b_t(\cla_0)$ where for each $t \in \IR$ we set $\cla^b_{t]}$ for 
the von-Neumann algebra $\{ j^b_s(x): s \le t, x \in \cla_0 \}''$. 

\vsp
\NI {\bf THEOREM 3.4: } Let $(\cla_0,\tau_t,\phi_0)$ be a Markov semigroup with a faithful normal invariant state $\phi_0$. Then the
following are equivalent:

\NI (a) $\phi_0(\tilde{\tau}_t(x)\tilde{\tau}_t(y)) \raro \phi_0(x)\phi_0(y)$ as $t \raro \infty$ for any $x,y \in \cla_0$.

\NI (b) $||\phi \tau_t - \phi_0|| \raro 0$ as $t \raro \infty$ for any normal state $\phi$ on $\cla_0$.

\NI {\bf PROOF:} For each $t \in \!R$ let $\cla^b_{t]}$ be the von-Neumann algebra generated by the backward processes
$\{ j_s^b(x): - \infty <  s \le t \}$ [Mo1]. Assume (a).  By Theorem 3.9 and Theorem 4.1 in [Mo1] we verify that weak$^*$ 
closure of $\bigcup_{t \in \!R } \cla^b_{t]}$ is $\clb(\clh)$. Since for each $t \in \!R$ the commutant of $\cla^b_{t]}$ contains
$\cla_{[t}$ we conclude that $\bigcap_{t \in \!R} \cla_{[t}$ is trivial. Hence (b)
follows once we appeal to Theorem 3.3. For the converse, it is enough if we verify that $\phi_0(\tilde{\tau}_t(x)J\tilde{\tau}_t(y)J) \raro
\phi_0(x)\phi_0(y)$ as $t \raro \infty$ for any $x,y \in \cla_0$ with $y \ge 0$ and $\phi_0(y)=1$.  To that end we check the following
easy steps $\phi_0(\tilde{\tau}_t(x)J\tilde{\tau}_t(y)J) = \phi_0(\tau_t(\tilde{\tau}_t(x))JyJ)$ and for any normal state $\phi$,
$|\phi \circ \tau_t(\tilde{\tau}_t(x)) - \phi_0(x)| \le ||\phi \circ \tau_t - \phi_0|| ||\tilde{\tau}_t(x)|| \le ||\phi \circ \tau_t-\phi_0||
||x||$. Thus the result follows once we note that $\phi$ defined by $\phi(x)=\phi_0(xJyJ)$ is a normal state. \qed

\vsp
\NI {\bf THEOREM 3.5:}  Let $(\cla_0,\tau_t,\phi_0)$ be a Markov semigroup with a normal invariant state $\phi_0$. Consider 
the following statements:

\NI (a) $\phi_0(\tau_t(x)\tau_t(y)) \raro \phi_0(x)\phi_0(y)$ as $t \raro \infty$ for all $x,y \in \cla_0$.

\NI (b) the strong $\mbox{lim}_{t \raro - \infty}F_{t]} = |\Omega><\Omega|$.

\NI (c) $\cla = \clb(\clh)$ 

Then (a) and (b) are equivalent statements and in such a case (c) is also true. If $\phi_0$ is also faithful (c) is also 
equivalent to (a) ( and hence ( b)).  

\vsp
\NI {\bf PROOF: } That (a) and (b) are equivalent is nothing but a restatement of Theorem 3.9 in [Mo1]. That (b) implies (c) is
obvious since the projection $[\cla'\Omega]$, where $\cla'$ is the commutant of $\cla$, is the support of the vector state
in $\cla$. We will prove now (c) implies (a). In case $\cla=\clb(\clh)$, we have $\bigcap_{t \in \!R} \cla^b_{t]}= \!C$,
thus in particular $\bigcap_{t \le 0} \alpha_t(\cla^b_{0]})=\!C$. Hence by Theorem 3.3 applied for the time-reverse endomorphism 
we verify that $||\phi \tilde{\tau}_t-\phi_0|| \raro 0$ as $t \raro \infty$. Now (a) follows once we appeal to Theorem 3.4 
for the adjoint semigroups since $\tilde{\tilde{\tau}}_t=\tau_t$. \qed 

\vsp
\NI {\bf THEOREM 3.6: } Let $(\cla_0,\tau_t,\phi_0)$ be as in Theorem 3.1. Then the following hold:

\NI (a) If $(\cla_0,\tau_t,\phi_0)$ is mixing then $\alpha_t(X) \raro \phi(X)$ as $t \raro \infty$ for
all $X \in \clb$, where $\clb$ is the $C^*$ completion of the $*$ algebra generated by 
$\{j_t(x): t \in \IR, x \in \cla_0 \}$. 

\NI (b) If $(\cla_0,\tau_t,\phi_0)$ is mixing and $\cla$ is a type-I factor then $\cla=\clb(\clh)$.  

\vsp
\NI{\bf PROOF: } For (a) we refer to [AM, Mo1]. By our hypothesis $\cla$ is a type-I von-Neumann factor and thus 
there exists an irreducible representation $\pi$ of $\clb$ in a Hilbert space $\clh_{\pi}$ quasi equivalent to $\pi_{\phi}$. 
There exists a density matrix $\rho$ on $\clh_{\pi}$ such that $\phi(X)= tr(\pi(X)\rho)$ for all $X \in \clb$. Thus there 
exists a unitary representation $t \raro U_t$ on $\clh_{\pi}$ so that $$U_t\pi(X)U^*_t=\pi(\alpha_t(X))$$ 
for all $t \in \IR$ and $X \in \clb$. Since $\phi= \phi \alpha_t$ on $\clb$ we also have $U^*_t \rho U_t=\rho$. 
We claim that $\rho$ is a one dimensional projection. Suppose not and then there exists at least two characteristic 
unit vectors $f_1,f_2$ for $\rho$ so that $f_1,f_2$ are characteristic vector for unitary representation $U_t$. Hence
we have $<f_i,\pi(X)f_i>= <f_i,\pi(\alpha_t(X))f_i>$ for all $t \in \IR$ and $i=1,2$. By taking limit we conclude
by (a) that $<f_i,\pi(X)f_i>=\phi(X)<f_i,f_i>=\phi(X)$ for $i=1,2$ for all $X \in \clb$. This violets irreducibility of 
representation $\pi$. \qed

\vsp
\NI {\bf PROPOSITION 3.7: } Let $(\cla_0,\tau_t,\phi_0)$ be as in Theorem 3.5 with $\phi_0$ as faithful. Then the commutant of 
$\cla_{[t}$ is $\cla^b_{t]}$ for each $t \in \IR$. 

\vsp
\NI {\bf PROOF: } It is obvious that $\cla_{[0}$ is a subset of the commutant of $\cla^b_{0]}$. Note also that $F_{[0}$
is an element in $\cla^b_{0]}$ which commutes with all the elements in $\cla_{[0}$. As a first step note that it is good 
enough if we show that $F_{[0}(\cla^b_{0]})'F_{[0} = F_{[0}\cla_{[0}F_{[0}$. As for some $X \in (\cla^b_{0]})'$ and 
$Y \in \cla_{[0}$ if we have $XF_{[0} = F_{[0}XF_{[0}=F_{[0}YF_{[0}=YF_{[0}$ then we verify that $XZf=YZf$ where $f$ 
is any vector so that $F_{[0}f=f$ and $Z \in \cla^b_{0]}$ and thus as such vectors are total in $\clh$ we get $X=Y$ ). Thus 
all that we need to show that $F_{[0}(\cla^b_{0]})'F_{[0} \subseteq F_{[0}\cla_{[0}F_{[0}$ as inclusion in other direction
is obvious. We will explore in following the relation that 
$F_{0]}F_{[0}=F_{[0}F_{0]}=F_{\{ 0 \}}$ i.e. the projection on the fiber at $0$ repeatedly. 
A simple proof follows once we use explicit formulas for $F_{0]}$ and $F_{[0}$ given in [Mo1].  

\vsp
Now we aim to prove that $F_{[0}\cla_{[0}'F_{[0} \subseteq F_{[0}\cla^b_{0]}F_{[0}$. Let 
$X \in F_{[0}\cla_{[0}'F_{[0}$ and verify that $X\Omega=XF_{0]}\Omega = F_{0]}XF_{0]}\Omega =F_{\{0\}}XF_{\{0\}}\Omega \in 
[j^b_0(\cla_0)''\Omega]$. On the other-hand we note by Markov property of the backward process 
$(j^b_t)$ that $F_{[0}\cla^b_{0]}F_{[0}=j^b(\cla_0)''$. Thus there exists an element $Y \in \cla^b_{0]}$ so that $X\Omega=Y\Omega$. 
Hence $XZ\Omega=YZ\Omega$ for all $Z \in \cla_{[0}$ as $Z$ commutes with both $X$ and $Y$. Since $\{ Z\Omega: Z \in \cla_{[0} \}$ 
spans $F_{[0}$, we get the required inclusion. Since inclusion in the other direction is trivial as 
$F_{[0} \in \cla_{[0}'$ we conclude that $F_{[0}\cla_{[0}'F_{[0} = F_{[0}\cla^b_{0]}F_{[0}.$

\vsp 
$F_{[0}$ being a projection in $\cla^b_{0]}$ we verify that $F_{[0}(\cla^b_{0]})'F_{[0} \subseteq (F_{[0}\cla^b_{0]}F_{[0})'$ 
and so we also have $F_{[0}(\cla^b_{0]})'F_{[0} \subseteq ( F_{[0} \cla_{[0}'F_{[0} )' $ as $\cla^b_{0]} \subseteq \cla_{[0}'$. 
Thus it is enough if we prove that  
$$F_{[0} \cla_{[0}'F_{[0} = (F_{[0} \cla_{[0} F_{[0})'$$
We will verify the non-trivial inclusion for the above equality. Let $X \in (F_{[0} \cla_{[0} F_{[0})'$ 
then $X\Omega=XF_{0]}\Omega=F_{0]}XF_{0]}\Omega = F_{\{0\}}XF_{\{0\}}\Omega  \in [j^b_0(\cla_0)\Omega]$. 
Hence there exists an element $Y \in F_{[0} \cla_{[0}'F_{[0}$ so that $X\Omega=Y\Omega$. Thus for any $Z \in \cla_{[0}$ 
we have $XZ\Omega=YZ\Omega$ and thus $XF_{[0}=YF_{[0}$. Hence $X=Y \in F_{[0} \cla_{[0}'F_{[0}$. Thus we get the required inclusion. 
  
\vsp
Now for any value of $t \in \IR$ we recall that $\alpha_t(\cla_{[0})=\cla_{[t}$ and $\alpha_t(\cla_{[0})'= 
\alpha_t(\cla_{[0}')$, $\alpha_t$ being an automorphism. This completes the proof as $\alpha_t(\cla^b_{0]})=\cla^b_{t]}$ 
by our construction.  \qed

\vsp
One interesting problem that we raised in [Mo1] whether Kolmogorov's property is time reversible i.e. whether 
$F_{t]} \raro |\Omega><\Omega|$ as $t \raro -\infty$ if and only if $F_{[t} \raro |\Omega><\Omega|$ as 
$t \raro \infty$. That it is true in classical case follows by Kolmogorov-Sinai-Rohlin theory on dynamical entropy 
for the associated Markov shift [Pa]. In the present general set up, it is true if $\cla_0$ is a type-I von-Neumann 
algebra with centre atomic [Mo1]. It is obviously true if the Markov semigroup is KMS symmetric. But in general it is 
false. In the last section we will give a class of counter example. This indicates that the quantum counter part of 
Kolmogorov property is unlikely to be captured by a suitable notion of quantum dynamical entropy with Kolmogorov-Sinai-Rohlin 
property.

\newsection{ Inductive limit state and purity: }

\vsp
Let $(\clb_0, \lambda_t,\;t \ge 0, \psi)$ be a unital $*-$ endomorphism with an invariant normal state $\psi$ on a 
von-Neumann algebra $\clb_0$ acting on a Hilbert space $\clh$. Let $P$ be the support projection for $\psi$. 
We set $\cla_0 = P \clb P$, a von-Neumann algebra acting on $\clh_0$, the closed subspace $P$, and $\tau_t(x)= P\lambda_t(PxP)P$, 
for any $x \in \cla_0$ and $t \ge 0$. Since $\lambda_t(P) \ge P$, it is simple to verify [Mo1] that 
$(\cla_0,\tau_t,\psi_0)$ is a quantum dynamical semigroup with a faithful normal invariant state $\psi_0$, 
where $\psi_0(x)=\psi(PxP)$ for $x \in \cla_0$. Now we set $j_0(x)=PxP$ and $j_t(x)=\lambda_t(j_0(x))$ for $t \ge 0$ and $x \in \cla_0$.  
A routine verification says that $F_{s]}j_t(x)F_{s]}=j_s(\tau_{t-s}(x))$ for $0 \le s \le t$, where $F_{s]}=\lambda_s(P),\; s \ge 0$. 
Let $\cla_{[0}$ be the von-Neumann algebra $\{ j_t(x): t \ge 0, x \in \cla_0 \}''$. As in Section 4 we check that 
$P\alpha_t(X)P=\tau_t(PXP)$ for all $X \in \cla_{[0}$. However are these vectors $\{ \lambda_{t_n}(PX_nP)....\lambda_{t_1}(PX_1P)f: \; 
f \in \clh_0, 0 \le t_1 \le t_2 \le ..\le t_k \le ..t_n, X_1,..,X_n \in \clb_0, n \ge 1 \}$ total in $\clh$? As an counter example in 
discrete time we consider an endomorphism on $\clb(\clh)$ [BJP] with a pure mixing state and note that $\cla_0$ is only scalers.  
Thus the cyclic space generated by the process $(j_t)$ on the pure state is itself. Thus the problem is rather delicate even when the 
von-Neumann algebra is the algebra of all bounded operators on $\clk$. We will not address this problem here. Since $\lambda_t(P) 
\lambda_{t_n}(PX_nP)...\lambda_{t_1}(PXP) \clh_0 = \lambda_{t_n}(PX_nP)...\lambda_{t_1}(PXP) \Omega $ for $t \ge t_n$, $\mbox{lim}_
{t \raro \infty}\lambda_t(P) =1$ is a necessary condition for cyclic property. The same counter example shows that it is not sufficient. 
In the following we explore the fact the support projection $P$ is indeed an element in the von-Neumann algebra $\cla$ generated by the 
process $(k_t(x):\;t \ge 0,\;x \in \cla_0)$ and asymptotic limit of the endomorphism $(\clb_0,\lambda_t, t \ge 0,\psi)$ is related 
with that of minimal endomorphism $(\cla_{[0},\alpha_t,\; t \ge 0 \phi)$. 

\vsp
In the following we consider a little more general situation. Let $\clb_0$ be a $C^*$ algebra, $(\lambda_t:\;t \ge 0)$ be a semigroup 
of injective endomorphisms and $\psi$ be an invariant state for $(\lambda_t:t \ge 0)$. We extend $(\lambda_t)$ to an automorphism on 
the $C^*$ algebra $\clb_{-\infty}$ of the inductive limit 
$$ \clb_0 \raro^{\lambda_t} \clb_0 \raro^{\lambda_t} \clb_0 $$
and extend also the state $\psi$ to $\clb_{-\infty}$ by requiring $(\lambda_t)$ invariance.
Thus there exists a directed set ( i.e. indexed by $\IT$ , by inclusion $\clb_{[-s} \subseteq \clb_{[-t}$ 
if and only if $t \ge s$ ) of C$^*$-subalgebras $\clb_{[t}$ of $\clb_{-\infty}$ so that the uniform closure of 
$\bigcup_{s \in \IT} \clb_{[s}$ is $\clb_{[-\infty}$. Moreover there exists an isomorphism
$$i_0: \clb_0 \raro \clb_{[0}$$
( we refer [Sa] for general facts on inductive limit of C$^*$-algebras). 
It is simple to note that $i_t=\lambda_t \circ i_0$ is an isomorphism of $\clb_0$ onto $\clb_{[t}$ and 
$$\psi_{-\infty} i_t = \psi$$ 
on $\clb_0$. Let $(\clh_{\pi},\pi,\Omega)$ be the GNS space associated with $(\clb_{[-\infty},\psi_{[-\infty})$ and $(\lambda_t)$ be 
the unique normal extension to $\pi(\clb_{-\infty})''$. Thus the vector state $\psi_{\Omega}(X)=<\Omega, X \Omega>$ is 
an invariant state for automorphism $(\lambda_t)$. As $\lambda_t(\clb_{[0}) \subseteq \clb_{[0}$ for all $t \ge 0$, 
$(\pi(\clb_{[0})'',\lambda_t,\;t \ge 0,\psi_{\Omega})$ is a quantum dynamics of endomorphisms. Let $F_{t]}$ be the support 
projection of the normal vector state $\Omega$ in the von-Neumann sub-algebra $\pi(\clb_{[t})''$. $F_{t]} \in \pi(\clb_{[t})'' 
\subseteq \pi(\clb_{[-\infty})''$ is a monotonically decreasing sequence of projections as $t \raro -\infty$. Let 
projection $Q$ be the limit. Thus $Q \ge [\pi(\clb_{[-\infty})'\Omega] \ge |\Omega><\Omega|$. So $Q=|\Omega><\Omega|$ ensures
that $\psi$ on $\clb_{[-\infty}$ is pure. We aim to investigate when $Q$ is pure i.e. $Q=|\Omega><\Omega|$. 

\vsp
To that end we set von-Neumann algebra $\cln_0=F_{0]}\pi(\clb_{[0})''F_{0]}$ and define family 
$\{ k_t: \cln_0 \raro \pi(\clb_{-\infty})'',\; t \in \IT \}$ of $*-$homomorphisms by
$$k_t(x) = \lambda_t(F_{0]}xF_{0]}),\;\; x \in \cln_0$$  
It is a routine work to check that $(k_t:t \in \IT)$ is the unique up to isomorphism ( in the
cyclic space of the vector $\Omega$ generated by the von-Neumann algebra $\{k_t(x): t \in \IT, x \in \cln_0 
\}$ ) forward weak Markov process associated with $(\cln_0,\eta_t,\psi_0)$ where 
$\eta_t(x)=F_{0]}\alpha_t(F_{0]}xF_{0]})F_{0]}$ for all $t \ge 0$. It is minimal once restricted to the cyclic space generated 
by the process. Thus $Q=|\Omega><\Omega|$ when restricted to the cyclic subspace of the process if and only if 
$\psi_0(\eta_t(x)\eta_t(y)) \raro \psi_0(x)\psi_0(y) $ as $t \raro \infty$ for all $x,y \in \cln_0$.

\vsp
\NI {\bf PROPOSITION 4.1: } Let $G_{0]}$ be the cyclic subspace of the vector $\Omega$ generated by $\pi(\clb_{[0})$. 

\NI (a) $G_{0]} \in \pi(\clb_{[0})'$ and the map $h: \pi(\clb_{[0})'' \raro G_{0]}\pi(\clb_{[0})''G_{0]}$ defined by
$X \raro G_{0]}XG_{0]}$ is an homomorphism and the range is isomorphic to $\pi_0(\clb_0)''$, where 
$(\clh_{\pi_0},\pi_0)$ is the GNS space associated with $(\clb_0,\psi)$.  

\NI (b) Identifying the range of $h$ with $\pi_0(\clb_0)''$ we have
$$ h \circ \lambda_t(X)=\lambda_t(h(X)) $$
for all $X \in \pi(\clb_{[0})''$ and $t \ge 0$. 

\NI (c) Let $P$ be the support projection of the state $\psi$ in von-Neumann algebra $\pi_0(\clb_0)''$ and 
$\cla_0=P\pi_0(\clb_0)''P$. We set $\tau_t(x)=P\lambda_t(PxP)P$ for all $t \ge 0,\; x \in \cla_0$ and 
$\psi_0(x)=\psi(PxP)$ for $x \in \cla_0$. Then 

\NI (i) $h(F_{0]})=P$ and $h(\cln_0)=\cla_0$;

\NI (ii) $h(\eta_t(x)) = \tau_t(h(x))$ for all $t \ge 0$.   

\vsp
\NI {\bf PROOF: } The map $\pi(X)\Omega \raro \pi_0(X)\Omega_0$ has an unitary extension which intertwines 
the GNS representation $(\clh_0,\pi_0)$ with the sub-representation of $\clb_{[0}$ on the cyclic 
subspace $G_{0]}$. Thus (a) follows. (b) is a simple consequence as $i_0: \clb_0 \raro \clb_{[0}$ is a 
$C^*$ isomorphism which covariant with respect to $(\lambda_t)$ for all $t \ge 0$ i.e $\lambda_t i_0(x)=i_0(\lambda_t(x))$
for all $x \in \clb_0$. That $h(F_{0]})=P$ is simple as $h$ is an isomorphism and thus also a normal map taking 
support projection $F_{0]}$ of the state $\psi$ in $\pi(\clb_{[0})''$ to support projection $P$ of the state $\psi$
in $\pi_0(\clb_0)''$. Now by homomorphism property of the map $h$ and commuting property with 
$(\lambda_t)$ we also check that $h(\cln_0)=h(F_{0]}\pi(\clb_{[0})''F_{0]})=P\pi_0(\clb_0)''P=\cla_0$ and 
$$h(\eta_t(x))=h(F_{0]})\lambda_t(h(F_{0]})h(x)h(F_{0]}))$$
$$=P\lambda_t(Ph(x)P)P=\tau_t(h(x))$$ for all $t \ge 0.$   

\vsp
\NI {\bf THEOREM 4.2: } $Q$ is pure if and only if $\phi_0(\tau_t(x)\tau_t(y)) \raro \phi_0(x)\psi_0(y)$ 
as $t \raro \infty$ for all $x,y \in \cla_0$. 

\vsp
\NI {\bf PROOF: } For any fix $t \in \IT$ since $k_t(\cla_0)= F_{t]}\pi(\clb_{[t})''F_{t]}$, for any $X \in \clb_{[t}$ we have 
$QX\Omega=QF_{t]}XF_{t]}\Omega=Qk_t(x)\Omega$ for some $x \in \cla_0$. Hence $Q=|\Omega><\Omega|$ if and only if $Q=|\Omega><\Omega|$ 
on the cyclic subspace generated by $\{ k_t(x),\;t \in \IT, x \in \cla_0 \}$. Theorem 3.5 says now that $Q=|\Omega><\Omega|$ if and
only if $\psi_0(\eta_t(x)\eta_t(y)) \raro \psi_0(x)\psi_0(y)$ as $t \raro \infty$ for all $x \in \cln_0$, Since $h$ is an
homomorphism and $h \eta_t(x)= \tau_t(h(x))$, we also have $h(\eta_t(x))\eta_t(y))=\tau_t(h(x))\tau_t(h(x))$. Since 
$\phi_0 \circ h = \psi_0$ we complete the proof.  \qed

\vsp
\NI {\bf COROLLARY 4.3: } $\psi_{[-\infty}$ is a pure state if $\phi_0(\tau_t(x)\tau_t(y)) \raro \phi_0(x)\psi_0(y)$
as $t \raro \infty$ for all $x,y \in \cla_0$.

\vsp
\NI {\bf PROOF: } It follows by Theorem 4.2 as $Q \le [\pi(\clb_{[-\infty})'\Omega] \le |\Omega><\Omega|$. \qed 

\vsp
Our analysis above put very little light whether the sufficient condition given in Corollary 4.3 is also necessary
for purity. We will get to this point in next section where we will deal with a class of examples.

\newsection{ Kolmogorov's property and pure translation invariant states: }

\vsp
Let $\omega$ be a translation invariant state on UHF$_d$ algebra $\cla=\otimes_{\IZ}M_d$ and $\omega'$ be the
restriction of $\omega$ to UHF$_d$ algebra $\clb_0=\otimes_{\IN} M_d$. There is a one to one correspondence between 
a translation invariant state $\omega$ and $\lambda$ (one sided shift ) invariant state $\omega'$ on UHF$_d$ algebra 
$\otimes_{\IN} M_d$. Powers's [Po] criteria easily yields that $\omega$ is a factor state if and only if $\omega'$ 
is a factor state. A question that comes naturally here which property of $\omega'$ is related with the purity of $\omega$. 
A systematic account of this question was initiated in [BJKW] inspired by initial success of [FNW1,FNW2,BJP] and a sufficient 
condition is obtained. In a recent article [Mo2] this line of investigation was further explored and we obtained a necessary and 
sufficient condition for a translation invariant lattice symmetric factor state to be pure and the criteria can be described 
in terms of Popescu elements canonically associated with Cuntz's representation. That the state is lattice symmetric played 
an important role in the duality argument used in the proof. 

\vsp
Here as an application of our general result, we aim now to find one more useful criteria for a translation invariant factor 
state $\omega$ on a one dimensional quantum spin chain $\otimes_{\IZ}M_d$ to be pure.  We also prove that purity of a lattice 
symmetric translation invariant state $\omega$ is equivalent to Kolmogorov's property of a Markov semigroup canonically 
associated with $\omega$.  

\vsp
First we recall that the Cuntz algebra $\clo_d ( d \in \{2,3,.., \} )$ is the universal $C^*$-algebra
generated by the elements $\{s_1,s_2,...,s_d \}$ subject to the relations:

$$s^*_is_j=\delta^i_j1$$
$$\sum_{1 \le i \le d } s_is^*_i=1.$$

\vsp
There is a canonical action of the group $U(d)$ of unitary $d \times d$ matrices on $\clo_d$ given by
$$\beta_g(s_i)=\sum_{1 \le j \le d}\overline{g^j_i}s_j$$
for $g=((g^i_j) \in U(d)$. In particular the gauge action is defined by
$$\beta_z(s_i)=zs_i,\;\;z \in \IT =S^1= \{z \in \IC: |z|=1 \}.$$
If UHF$_d$ is the fixed point subalgebra under the gauge action, then UHF$_d$ is the closure of the
linear span of all wick ordered monomials of the form
$$s_{i_1}...s_{i_k}s^*_{j_k}...s^*_{j_1}$$
which is also isomorphic to the UHF$_d$ algebra
$$M_{d^\infty}=\otimes^{\infty}_1M_d$$
so that the isomorphism carries the wick ordered monomial above into the matrix element
$$e^{i_1}_{j_1}(1)\otimes e^{i_2}_{j_2}(2) \otimes....\otimes e^{i_k}_{j_k}(k) \otimes 1 \otimes 1 ....$$
and the restriction of $\beta_g$ to $UHF_d$ is then carried into action
$$Ad(g)\otimes Ad(g) \otimes Ad(g) \otimes ....$$

\vsp
We also define the canonical endomorphism $\lambda$ on $\clo_d$ by
$$\lambda(x)=\sum_{1 \le i \le d}s_ixs^*_i$$
and the isomorphism carries $\lambda$ restricted to UHF$_d$ into the one-sided shift
$$y_1 \otimes y_2 \otimes ... \raro 1 \otimes y_1 \otimes y_2 ....$$
on $\otimes^{\infty}_1 M_d$. Note that $\lambda \beta_g = \beta_g \lambda $ on UHF$_d$.

\vsp
Let $d \in \{2,3,..,,..\}$ and $\IZ_d$ be a set of $d$ elements.  $\cli$ be the set of finite sequences
$I=(i_1,i_2,...,i_m)$ where
$i_k \in \IZ_d$ and $m \ge 1$. We also include empty set $\emptyset \in \cli$ and set $s_{\emptyset }=1=s^*_{\emptyset}$,
$s_{I}=s_{i_1}......s_{i_m} \in \clo_d $ and $s^*_{I}=s^*_{i_m}...s^*_{i_1} \in \clo_d$.

\vsp
Let $\omega$ be a translation invariant state on $\cla= \otimes_{\IZ} M_d$ where $M_d$ is $(d \times d)$ matrices with 
complex entries. Identifying $\otimes_{\IN}M_d$ with $\mbox{UHF}_d$ we find a one to one relation from a $\lambda$ 
invariant state on $\mbox{UHF}_d$ with that of an one sided shift invariant state on $\cla_R=\otimes_{\IN} M_d$.    
Let $\omega'$ be an $\lambda$-invariant state on the $\mbox{UHF}_d$ sub-algebra of $\clo_d$. Following 
[BJKW, section 7], we consider the set
$$K_{\omega'}= \{ \psi: \psi \mbox{ is a state on } \clo_d \mbox{ such that } \psi \lambda =
\psi \mbox{ and } \psi_{|\mbox{UHF}_d} = \omega' \}$$
By taking invariant mean on an extension of $\omega'$ to $\clo_d$, we verify that $K_{\omega'}$ is non empty and
$K_{\omega'}$ is clearly convex and compact in the weak topology. In case $\omega'$ is an ergodic state ( extremal state )
$K_{\omega'}$ is a face in the $\lambda$ invariant states. Before we recall Proposition 7.4 of [BJKW] in the following
proposition.

\vsp
\NI {\bf PROPOSITION 5.1:} Let $\omega'$ be ergodic. Then $\psi \in K_{\omega'}$ is an extremal point in
$K_{\omega'}$ if and only if $\hat{\omega}$ is a factor state and moreover any other extremal point in $K_{\omega'}$
have the form $\psi \beta_z$ for some $z \in \IT$.

\vsp
We fix any $\hat{\omega} \in K_{\omega'}$ point and consider the associated Popescu system $(\clk,\clm,v_k,\Omega)$
described as in Proposition 2.4. A simple application of Theorem 3.6 in [Mo2] says that the inductive limit state
$\hat{\omega}_{-\infty}$ on the inductive limit $(\clo_d,\hat{\omega}) \raro^{\lambda} (\clo_d,\hat{\omega})
\raro^{\lambda} (\clo_d,\hat{\omega})$ is pure if $\phi_0(\tau_n(x)\tau_n(y)) \raro \phi_0(x)\phi_0(y)$
for all $x,y \in \clm$ as $n \raro \infty$. This criteria is of limited use in determining purity of $\omega$ unless
we have $\pi_{\hat{\omega}}(\mbox{UHF}_d)''=\pi_{\hat{\omega}}(\clo_d)''$. We prove a more powerful criteria in
the next section, complementing a necessary and sufficient condition obtained by [Mo2], for a translation 
invariant factor state $\omega$ to be pure.  

\vsp
To that end note that the von-Neumann algebra $\{ S_IS^*_J: |I|=|J| < \infty \}''$ acts on the cyclic subspace
of $\clh_{\pi_{\hat{\omega}}}$ generated by the vector $\Omega$. This is isomorphic with the GNS representation
associated with $(\clb_0,\omega')$. The inductive limit $(\clb_{-\infty},\hat{\omega}_{-\infty})$ [Sa] described as
in Proposition 3.6 in [Mo2] associated with $(\clb_0,\lambda_n,\;n \ge 0,\omega')$ is UHF$_d$ algebra
$\otimes_{\IZ}M_d$ and the inductive limit state is $\omega$.
\vsp
Let $Q$ be the support projection of the state $\hat{\omega}$ in $\pi_0(\clb_0)''$ and $\cla_0=Q\pi(\clb_0)''Q$. Since
$\psi_{\Omega}(\Lambda(X))=\psi_{\Omega}(X)$ for all $X \in \pi_{\hat{\omega}}(\mbox{UHF}_d)''$, $\Lambda(Q) \in
\pi_{\hat{\omega}}(\mbox{UHF}_d)''$ and $\Lambda(Q) \ge Q$ [Mo1]. Thus $Q\Lambda(I-Q)Q=0$ and we have
$(I-Q)S_k^*Q=0$ for all $1 \le k \le d$. The reduced Markov map $\eta:\cla_0 \raro \cla_0$ is defined by
\be
\eta(x)=Q\Lambda(QxQ)Q
\ee
for all $x \in \cla_0$ which admits a faithful normal state $\phi_0$ defined by
\be
\psi_0(x)=\psi_{\Omega}(QxQ),\;\; x \in \cla_0
\ee
\vsp
In particular, $\Lambda_n(Q) \uparrow I$ as $n \raro \infty$. Hence $ \{ S_If: |I| < \infty,\; Qf=f,f \in \clh_{\pi}
\}$
is total in $\clh_{\pi_{\hat{\omega}}}$.

\vsp
We set $l_k=QS_kQ$, where $l_k$ need not be an element in $\cla_0$. However $l_Il^*_J \in \cla_0$ provided $|I|=|J| < \infty$.
Nevertheless we have $Q\Omega=\Omega$ and thus verify that
$$\hat{\omega}(s_Is^*_J)=<\Omega,S_IS^*_J\Omega>$$
$$<\Omega,QS_IS^*_JQ\Omega>=<\Omega,l_Il^*_J\Omega>$$
for all $|I|,|J| < \infty$. In particular we have
$$\omega'(s_Is^*_J)=\psi_0(l_Il^*_J)$$
for all $|I|=|J| < \infty$.
\vsp
For each $n \ge 1$ we note that $\{ S_IS^*_J: |I|=|J| \le n \}'' \subseteq \Lambda_n(\pi_{\hat{\omega}}(\mbox{UHF}_d)'')'
\bigcap  \pi_{\hat{\omega}}(\mbox{UHF}_d)''$ and thus $\pi_{\hat{\omega}}(\mbox{UHF}_d)'' \subseteq
(\bigcap_{n \ge 1}\Lambda_n(\pi_{\hat{\omega}}(\mbox{UHF}_d)'')'$.
Hence
\be
\bigcap_{n \ge 1}\Lambda_n(\pi_{\hat{\omega}}(\mbox{UHF}_d)'') \subseteq
\pi_{\hat{\omega}}(\mbox{UHF}_d)'' \bigcap \pi_{\hat{\omega}}(\mbox{UHF}_d)'.
\ee

\vsp
Now by Proposition 1.1 in [Ar, see also Mo2] $||\psi \Lambda^n - \psi_{\Omega}|| \raro 0$
as $n \raro \infty$ for any normal state $\psi$ on $\pi_{\hat{\omega}}(\mbox{UHF}_d)''$
if $\omega'$ is a factor state. Thus we have arrived at the
following well-known result of R. T. Powers [Pow1,BR].
\vsp
\NI {\bf THEOREM 5.2: } Let $\omega'$ be a $\lambda$ invariant state on $\mbox{UHF}_d$ $\otimes_{\IN}M_d$. Then
the following statements are equivalent:

\vsp
\NI (a) $\omega'$ is a factor state;

\NI (b) For any normal state $\psi$ on $\cla_0$, $||\psi \eta_n-\psi_0|| \raro 0$ as $n \raro \infty;$

\NI (c) For any $x \in  \mbox{UHF}_d$ $\otimes_{\IN}M_d$
$$\mbox{sup}_{||y|| \le 1}|\omega'(x\lambda_n(y))-\omega'(x)\omega'(y)| \raro 0$$
as $n \raro \infty$;

\NI (d) $\omega'(x \lambda_n(y)) \raro \omega'(x)\omega'(y)$ as $n \raro \infty$ for all $x,y \in \mbox{UHF}_d$
$\otimes_{\IN}M_d$;

\vsp
\NI {\bf PROOF: } For any normal state $\psi$ on $\cla_0$ we note that $\psi_P(X)=\psi(PXP)$ is a normal
state on $\pi_{\hat{\omega}}(\mbox{UHF}_d)''$ and $||\psi \eta_n-\psi_0|| \le
||\psi_P \Lambda_n-\psi_{\Omega}||$. Thus by the above argument (a) implies (b). That (c) implies (d) and
(d) implies (a) are obvious. We will prove that (b) implies (c).
Note that for (c) it is good enough if we verify for all non-negative
$x \in \mbox{UHF}_d$ with finite support and $\omega'(x)=1$. In such a case for large values of $n$ the
map $\pi_{\hat{\omega}}(y) \raro \omega'(x \lambda_n(y))$ determines a normal state on $\pi_{\hat{\omega}}(\mbox{UHF}_d)''$.
Hence (c) follows whenever (b) hold. \qed

\vsp
\NI {\bf COROLLARY 5.3: } Let $\omega$ be a translation invariant state on $\mbox{UHF}_d$ $\otimes_{\IZ}M_d$. Then the
following are equivalent:

\vsp
\NI (a) $\omega$ is a factor state;

\NI (b) $\omega(x \lambda_n(y)) \raro \omega(x)\omega(y)$ as $n \raro \infty$ for all $x,y \in \mbox{UHF}_d$ $\otimes_{\IZ}M_d$;

\vsp
\NI {\bf PROOF: } First we recall $\omega$ is a factor state if and only if $\omega$ is an extremal point in the
translation invariant state i.e. $\omega$ is an ergodic state for the translation map. Since the cluster property
(b) implies ergodicity, (a) follows. For the converse note that $\omega$ is a ergodic state for the translation
map if and only if $\omega'$ is ergodic for $\lambda$ on $\mbox{UHF}_d$ $\otimes_{\IN}M_d$. Hence by Theorem 3.2
we conclude that statement (b) hold for any local elements $x,y \in \mbox{UHF}_d$ $\otimes_{\IZ}M_d$. Now we use
the fact that local elements are dense in the $C^*$ norm to complete the proof. \qed

\vsp
\NI {\bf PROPOSITION 5.4:} Let $\omega$ be a translation invariant extremal state on $\cla$ and $\psi$ be an extremal 
point in $K_{\omega}$. Then following hold:

\NI (a) $H =\{ z \in S^1: \psi \beta_z =\psi \}$ is a closed subgroup of $S^1$ and $\pi(\clo_d)''^{\beta_H}=\pi(\mbox{UHF}_d)''$. 
Furthermore we have $\bigcap_{n \ge 1} \Lambda^n(\pi(\clo_d)'') = \pi(\clo_d)'' \bigcap \pi(\mbox{UHF}_d)'$;

\NI (b) If $H=S^1$ then $\pi(\clo_d)'' \bigcap \pi(\mbox{UHF}_d)' = \IC$; 

\NI (c) Let $(\clh,\pi,\Omega)$ be the GNS representation of $(\clo_d,\psi)$ and $P$ be the support projection of the state 
$\psi$ in $\pi(\clo_d)''$. Then $P \in \pi(\mbox{UHF}_d)''$ is also the support projection of the state $\psi$ in 
$\pi(\mbox{UHF}_d)''$;

\vsp
\NI {\bf PROOF: } First part of (a) is noting but a restatement of Proposition 2.5 in [Mo2] modulo the factor property
of $\pi(\mbox{UHF}_d)''$. For a proof of the factor property we refer to Lemma 7.11 in [BJKW] modulo a modification
described in Proposition 3.2 in [Mo2].

\vsp
We aim now to show that $\bigcap_{n \ge 1} \Lambda^n(\pi(\clo_d)'') = \pi(\clo_d)'' \bigcap \pi(\mbox{UHF}_d)'$. It is obvious by Cuntz
relation that $\bigcap_{n \ge 1} \Lambda^n(\pi(\clo_d)'') \subseteq \pi(\clo_d)'' \bigcap \pi(\mbox{UHF}_d)'$. For the converse
let $X \in \pi(\clo_d)'' \bigcap \pi(\mbox{UHF}_d)'$ and fix any $n \ge 1$ and set $Y_n= S^*_IXS_I$ with $|I|=n$. Since $X
\in \pi(\mbox{UHF}_d)'$ we verify that $S^*_IXS_I=S^*_IXS_IS^*_JS_J=S^*_IS_IS^*_JXS_J=S^*_JXS_J$ for any $|J|=n$. Thus $Y_n$
is independent of the multi-index that we choose. Once gain as $X \in \pi(\mbox{UHF}_d)'$ we also check that
$\Lambda^n(Y_n)= \sum_{J:|J|=n} S_JS_I^*XS_IS_J^* = X$. Hence $X \in \bigcap_{n \ge 1} \Lambda^n(\pi(\clo_d)'')$.

\vsp
Now $\pi(\mbox{UHF}_d)''$ being a factor, a general result in [BJKW, Lemma 7.12] says that
$\pi(\clo_d)'' \bigcap \pi(\clo_d^H)'$ is a commutative von-Neumann algebra generated by
an unitary operator $u$ so that $\beta_z(u)=\gamma(z)u$ for all $z \in H$ and some character $\gamma$ of $H$.
Furthermore there exists a $z_0 \in H$ so that $\beta_{z_0}(x)=uxu^*$ for all $x \in \pi(\clo_d)''$.
Thus we also have $\beta_{z_0}(u)=u=\gamma(z_0)u$. So we have $\gamma(z_0)=1$. $H$ being $S^1$ the character
can be written as $\gamma(z)=z^k$ all $z \in H$ and for some $k \ge 1$. Hence $u^kx(u^k)^*= \beta_{z_0^k}(x)=x$.
$\pi(\clo_d)''$ being a factor $u^k$ is a scaler. By multiplying a proper factor we can choose an unitary
$u \in \pi(\clo_d)''\bigcap \pi(\mbox{UHF}_d)'$ so that $u^k=1$. However we also check that for all
$z \in S^1$ we have $\beta_z(u^k)=\gamma(z)^ku^k$ i.e. $\gamma(z)^k=1$ for all $z \in S^1$ as $u^k=1$. Hence
$\gamma(z)=1$ for all $z \in S^1$. Thus $\beta_z(u)=u$ for all $z \in S^1$ and $u$ is scaler as $u$ is also an
element in $\pi(\mbox{UHF}_d)''$ by the first part. $\pi(\mbox{UHF}_d)''$ being a factor we conclude that $u$
is a scaler. Hence $\pi_{\psi}(\clo_d)'' \bigcap \pi(\mbox{UHF}_d)'$ is trivial. This completes the proof of (b).

\vsp
It is obvious that $\beta_z(P)=P$ for all $z \in H$ and thus by (a) $P \in \pi(\mbox{UHF}_d)''$ and thus also the support
projection in $\pi(\mbox{UHF}_d)''$ of the state $\psi$. (c) is a simple consequence of (a) and Corollary 4.3.  \qed

\vsp
\NI {\bf THEOREM 5.4: } Let $\omega$ be a translation invariant state on $\mbox{UHF}_d$ $\otimes_{\IZ}M_d$ and $P$ be the
support projection of $\psi \in K_{\omega'}$ in $\pi(\clo_d)''$. Further let $\cla_0$ be the von-Neumann algebra $P\pi(\mbox{UHF}_d)''P$ 
acting on the subspace $P$ and completely positive map $\tau: \cla_0 \raro \cla_0$ defined by $\tau(x)=P\Lambda(PxP)P$, i.e.
$\tau(x)= \sum_k l_k x l_k^*$ be the completely positive map on $\cla_0$ where $l_k=P\pi(s_k)P$ for all $1 \le k \le d$. Then the 
following hold:

\NI (a) If $\phi_0(\tau^n(x)\tau^n(y)) \raro \phi_0(x)\phi_0(y)$ as $n \raro \infty$ for all $x,y \in \cla_0$ then $\omega$ is pure;

\NI (b) If $H=S^1$ then $||\phi\tau^n - \phi_0||  \raro 0$ as $n \raro \infty$ for any normal state on $\cla_0$;

\vsp
\NI {\bf PROOF: } (a) follows by an easy application of Corollary 4.3. For a proof for (b) we appeal to [Ar, Proposition 1.1] 
and the last statement in Proposition 5.3 (a). \qed

\vsp
By a duality argument, Theorem 3.4 in [Mo2], $||\psi \eta_n-\psi_0|| \raro 0$ as $n \raro
\infty$ for any normal state $\psi$ if and only if $|\psi_0(\tilde{\eta}_n(x)\tilde{\eta}_n(y))
\raro \psi_0(x)\psi_0(y)|$ as $n \raro \infty$ for any $x,y \in \cla_0$, where $(\cla_0,\tilde{\eta},\phi_0)$
the KMS-adjoint Markov semigroup [OP,AcM,Mo1] of $(\cla_0,\eta,\phi_0)$. 

\vsp
We recall the unique KMS state $\psi=\psi_{\beta}$ on $\clo_d$ where $\beta=ln(d)$ is a factor state and 
$\psi_{\beta} \in K_{\omega}$ where $\omega'$ is the unique trace on $\mbox{UHF}_d$. For a proof that $H=S^1$ for $\psi_{\beta}$ 
we refer to [BR]. $\omega$ is the unique trace on $\cla$ and so is a factor state. Hence by Proposition 5.4 (d) we have 
$\pi_{\psi}(\clo_d)'' \bigcap \pi_{\psi}(\mbox{UHF}_d)'$ is trivial. Thus $\bigcap_{n \ge 1 } \Lambda(\pi_{\psi}
(\clo_d)'') = \IC$. In particular $\bigcap_{n \ge 1} \Lambda^n(\pi_{\psi}(\mbox{UHF}_d)'') = \IC$. On the other hand 
$\psi_{\beta}$ being faithful, the support projection is the identity operator and thus canonical Markov semigroup $\tau$ is 
equal to $\Lambda$. $\Lambda$ being an endomorphism and $\psi_{\beta}$ being faithful, we easily verify 
that $\tau$ does not admit Kolmogorov property. On the other hand $H=S^1$ and so by Proposition 5.4 (d)  
$||\phi \tau^n - \phi_0|| \raro 0$ as $n \raro \infty$ for any normal state $\phi$ on $\cla_0$. This
example unlike in the classical case shows that Kolmogorov's property of a non-commutative 
dynamical system in general is not time reversible.

\bigskip
{\centerline {\bf REFERENCES}}

\begin{itemize} 

\bigskip
\item{[AM]} Accardi, L., Mohari, A.: Time reflected Markov processes. Infin. Dimens. Anal. Quantum Probab. Relat. Top., vol-2, no-3, 397-425 (1999).

\item {[Ar]} Arveson, W.: Pure $E_0$-semigroups and absorbing states, Comm. Math. Phys. 187 , no.1, 19-43, (1997)

\item {[BP]} Bhat, B.V.R., Parthasarathy, K.R.: Kolmogorov's existence theorem for Markov processes on $C^*$-algebras, Proc.
Indian Acad. Sci. 104,1994, p-253-262.

\item {[BR]} Bratteli, Ola., Robinson, D.W. : Operator algebras and quantum statistical mechanics, I,II, Springer 1981.

\item {[BJ]} Bratteli, Ola; Jorgensen, Palle E. T. Endomorphism of $\clb(\clh)$, II,
Finitely correlated states on $\clo_N$, J. Functional Analysis 145, 323-373 (1997).

\item {[BJP]} Bratteli, Ola., Jorgensen, Palle E.T. and Price, G.L.: Endomorphism of $\clb(\clh)$, Quantization, nonlinear partial differential 
equations, Operator algebras, ( Cambridge, MA, 1994), 93-138, Proc. Sympos. Pure Math 59, Amer. Math. Soc. Providence, RT 1996.

\item {[BJKW]} Bratteli, O., Jorgensen, Palle E.T., Kishimoto, Akitaka and
Werner Reinhard F.: Pure states on $\clo_d$, J.Operator Theory 43 (2000),
no-1, 97-143.

\item{[Da]} Davies, E.B.: Quantum Theory of open systems, Academic press, 1976.

\item{[FNW1]} Fannes, M., Nachtergaele,D., Werner,R.: Finitely Correlated States on Quantum Spin Chains, Commun. Math. Phys. 144,
443-490 (1992).

\item{[FNW2]} Fannes, M., Nachtergaele,D., Werner,R.: Finitely Correlated pure states, J. Funct. Anal. 120, 511-534
(1994).

\item{[Fr]} Frigerio, A.: Stationary states of quantum dynamical
semigroups. Commun. Math. Phys. 63, 269-276 (1978).

\item{[Li]} Lindblad, G. : On the generators of quantum
dynamical semigroups, Commun.  Math. Phys. 48, 119-130 (1976).

\item{[Mo1]} Mohari, A.: Markov shift in non-commutative probability, Jour. Func. Anal. 199 (2003) 189-209.  

\item{[Mo2]} Mohari, A.: $SU(2)$ symmetry breaking in quantum spin chain, The preprint is under review in 
Communication in Mathematical Physics, http://arxiv.org/abs/math-ph/0509049.
 
\item{[Mo3]} Mohari, A.: Quantum detailed balance and split property in quantum spin chain, Arxiv:  
http://arxiv.org/abs/math-ph/0505035.

\item{ [Mo4]} Mohari, A: Jones index of a Markov semigroup, Preprint 2007.

\item{[Mo5]} Mohari, A.: Ergodicity of Homogeneous Brownian flows, Stochastic Process. Appl. 105 (1),99-116.

\item{[OP]} Ohya, M., Petz, D.: Quantum entropy and its use, Text and monograph in physics, Springer-Verlag 1995. 

\item{[Po]} Powers, Robert T.: An index theory for semigroups of $*$-endomorphisms of
$\clb(\clh)$ and type II$_1$ factors.  Canad. J. Math. 40 (1988), no. 1, 86--114.

\item{[Pa]} Parry, W.: Topics in Ergodic Theory, Cambridge University Press, Cambridge, 1981. 

\item{[Sak]} Sakai, S.: C$^*$-algebras and W$^*$-algebras, Springer 1971.  

\item{[Sa]} Sauvageot, Jean-Luc: Markov quantum semigroups admit covariant Markov $C^*$-dilations. Comm. Math. Phys. 
106 (1986), no. 1, 91103.

\item{[So]} Stormer, Erling : On projection maps of von Neumann algebras. Math. Scand. 30 (1972), 46--50.

\item{[Vi]} Vincent-Smith, G. F.: Dilation of a dissipative quantum dynamical system to a quantum Markov process. Proc. 
London Math. Soc. (3) 49 (1984), no. 1, 5872. 

\end{itemize}

\end{document}